\documentclass[11pt]{article}
\usepackage{psfrag,graphicx,
}
\usepackage{amsfonts,amssymb}
\usepackage{fullpage,amsmath,amscd,amsthm,amssymb,amsxtra,latexsym,epsf, 
epic,eepic,mathrsfs}
\usepackage{epsf,latexsym}
\input xy
\xyoption{all}
\input epsf
\DeclareMathOperator{\Der}{Der}
\DeclareMathOperator{\Ext}{Ext}
\DeclareMathOperator{\Gl}{Gl}
\DeclareMathOperator{\Sl}{Sl}
\DeclareMathOperator{\Rep}{Rep}
\DeclareMathOperator{\indeg}{\mathbf {in}}
\DeclareMathOperator{\outdeg}{\mathbf {out}}

\begin{document}
\newfont{\Bb}{msbm10}
\newcommand\s{\mathscr}

\def\cc{\circ}
\def\y{\item}
\def\ni{\noindent}
\def\ben{\begin{enumerate}}
\def\een{\end{enumerate}}
\def\beq{\begin{equation}}
\def\eeq{\end{equation}}
\def\bit{\begin{itemize}}
\def\eit{\end{itemize}}
\def\ec{\end{center}}
\def\bc{\begin{center}}
\def\ld{{\ldots}}
\def\cd{{\cdots}}


\def\d{\mathbf{d}}
\def\e{\mathbf{e}}
\def\L{\mathfrak{l}}
\def\pgl{\mathfrak{pgl}}
\def\gl{\mathfrak{gl}}
\def\gg{\mathfrak{g}}
\def\la{\langle}
\def\ra{\rangle}
\def\ds{\displaystyle}

\def\th{{\theta}}
\def\ve{{\varepsilon}}
\def\pp{{\varphi}}
\def\a{{\alpha}}
\def\l{{\lambda}}
\def\rr{{\rho}}
\def\ss{{\sigma}}
\def\g{{\gamma}}
\def\G{{\Gamma}}
\def\O{{\Omega}}
\def\k{{\kappa}}


\def\Proof{{\bf Proof}\hspace{.2in}}
\def\eop{\hspace*{\fill}$\Box$ \vskip \baselineskip}
\def\eoq{{\hfill{$\Box$}}}


\def\tr{{\pitchfork}}
\def\b{{\bullet}}
\def\w{{\wedge}}
\def\p{{\partial}}


\def\supp{\mbox{supp}}
\def\im{\mbox{im}}

\def\Mat{\mbox{Mat}}
\def\sgn{\mbox{sign}}
\def\cod{\mbox{codim }\ }
\def\hh{\mbox{height }\ }
\def\opp{\mbox{\scriptsize opp}}
\def\vol{\mbox{vol}}
\def\deg{\mbox{deg}}
\def\Sym{\mbox{Sym}}
\def\Hom{\mbox{Hom}}
\def\End{\mbox{End}}
\def\Aut{\mbox{Aut}}
\def\DD{{\Der(\log D)}}
\def\Dh{{\Der(\log h)}}
\def\pd{\mbox{pd}}
\def\dim{{\mbox{dim}}}
\def\depth{{\mbox{depth}}}


\def\AA{\mbox{\Bb A}}
\def\QQ{\mbox{\Bb Q}}
\def\CC{\mbox{\Bb C}}
\def\RR{\mbox{\Bb R}}
\def\HH{\mbox{\Bb H}}
\def\PP{\mbox{\Bb P}}
\def\ZZ{\mbox{\Bb Z}}
\def\NN{\mbox{\Bb N}}
\def\EE{{\cal E}}
\def\R{{\cal R}}
\def\C{{\cal C}}
\def\A{{\cal A}}
\def\TT{{\cal T}}
\def\OO{{\cal O}}
\newcommand{\bbbz}{{\mathbb Z}}


\newcommand{\n}[1]{\| #1 \|}
\newcommand{\um}[1]{{\underline{#1}}}
\newcommand{\om}[1]{{\overline{#1}}}
\newcommand{\fl}[1]{\lfloor #1 \rfloor}
\newcommand{\ce}[1]{\lceil #1 \rceil}
\newcommand{\ncm}[2]
{{\left(\!\!\!\begin{array}{c}#1\\#2\end{array}\!\!\!\right)}}
\newcommand{\ncmf}[2]
{{\left[\!\!\!\begin{array}{c}#1\\#2\end{array}\!\!\!\right]}}
\newcommand{\ncms}[2]
{{\left\{\!\!\!\begin{array}{c}#1\\#2\end{array}\!\!\!\right\}}}


\def\iff{{\Leftrightarrow}}
\def\imp{{\Rightarrow}}
\def\to{{\ \rightarrow\ }}
\def\too{{\ \longrightarrow\ }}
\def\into{{\hookrightarrow}}
\def\st{\stackrel}


\newtheorem{theorem}{Theorem}[section]
\newtheorem{lemma}[theorem]{Lemma}
\newtheorem{sit}[theorem]{}
\newtheorem{lemmadefinition}[theorem]{Lemma and Definition}
\newtheorem{proposition}[theorem]{Proposition}
\newtheorem{example}[theorem]{Example}
\newtheorem{corollary}[theorem]{Corollary}
\newtheorem{definition}[theorem]{Definition}
\newtheorem{conjecture}[theorem]{Conjecture}
\title{Linear Free Divisors and Quiver Representations}
\author{Ragnar-Olaf Buchweitz and David Mond}
\maketitle
\hfill To Gert-Martin Greuel\\
\section{Introduction}
Let $X$ be a non-singular $n$-dimensional complex
manifold (or 
algebraic variety over 
an algebraically closed field $k$ of characteristic zero),
and let $D\subset X$ be a hypersurface with reduced defining ideal
$I_X$.
We denote by $\Der(-\log D)$ the
sheaf
 of vector fields $\chi \in \Der_{X 
}$ such that $\chi\cdot I_X\subset I_X$, or, equivalently, 
such that $\chi$ is tangent to $D$ at its regular points. It is clearly
an $\OO_X$-module.\\
\begin{definition}
The hypersurface $D\subseteq X$ is a {\em free divisor\/} if $\Der(-\log D)$ 
is a locally free $\OO_X$-module.
\end{definition}
Free divisors were introduced by K.\,Saito in \cite{saito}.
The simplest example is the normal crossing divisor, but the main source
of examples, motivating Saito's definition, has been the deformation theory of 
singularities, where discriminants and bifurcation sets are 
 frequently
free divisors. If $D$ is the discriminant hypersurface in the base $S$ of a
versal deformation of an isolated hypersurface singularity, the module
$\Der(-\log D)$ is the kernel of the Kodaira-Spencer map from 
$\Der_{S}$
onto the relative $T^1$ of the deformation, and from this freeness follows by
an easy homological argument, due initially to Teissier. 
Variants on this argument show the freeness of the discriminant  
in the base of 
a versal deformation in a number of cases: isolated
complete interesection singularities (\cite{looijenga}), space-curve
singularities (\cite{duco}), functions on space curves (\cite{goryunov},
\cite{ms}), Gorenstein surface singularities in 5-space 
(\cite{buchweitzebeling}),
Hilbert schemes of a smooth surface (\cite{buchweitz}). Damon, in his paper 
``The legacy of free divisors'' (\cite{damonleg}), has shown, by
an essentially similar argument, how the 
bifurcation set in the base space of a versal deformation of a 
non-linear section of a free divisor is once again a free divisor, 
provided a natural condition, namely, the existence of 
``Morse-type singularities'', is met.
Another significant source of examples is the theory of hyperplane 
arrangements, where many examples of free arrangements have been
constructed by combinatorial means (see e.g. \cite{orlikterao} Chapter 4).

Saito's original paper \cite{saito} contained the following criterion, 
now known by his name, for a divisor $D$ to be free: 
\begin{proposition}{\em ({\sc Saito's Criterion})}
\label{scrit} 

The hypersurface $D\subset X$ is a free divisor in
the neighbourhood of a point $x$ if and 
only if there are germs  
of vector fields $\chi_1,\ld,\chi_n \in\Der(-\log D)_{x}$,
such that the determinant of the matrix of coefficients 
$[\chi_1,\ld,\chi_n]$,
with respect to some, or any, 
$\OO_{X,x}$-basis of $\Der_{X,x}$, is a reduced equation for $D$ at $x$. 
In this case, $\chi_1,\ld, \chi_n$ form a basis for $\Der(-\log D)_{x}$.\eop
\end{proposition}
Note that it is clear that the determinant of the 
matrix of coefficients of any $n$-tuple
of vector fields in $\Der(-\log D)$ must vanish identically on $D$, since 
at any regular point $x\in D$
all $n$ vectors lie in the $n-1$-dimensional vector space $T_xD$.
 Moreover,
since $\Der(-\log D)$ coincides with $\Der_X$ outside $D$, the 
determinant of the matrix of coefficients 
of any set of generators of $\Der(-\log D)$ must vanish
only on $D$.

In practice, one uses often the following concrete algebraic version of this 
criterion that does not refer to vector fields directly, rather characterizes 
the Taylor series of the function $f$ defining a free divisor at some point 
$x\in X$:

\begin{proposition}
A formal power series $f\in P=k[\![z_{1},...,z_{n}]\!]$ defines a (formal) 
free divisor, if it is reduced, that is, squarefree, and there is an 
$(n{\times}n)$--matrix $A$ with entries from $P$ such that
$$
\det A = f\quad\text{and}\quad 
(\nabla f)A\equiv (0,...,0) \bmod f\,,
$$
where $\nabla f =
\left(
\frac{\partial f}{\partial z_{1}},\ldots,\frac{\partial f}{\partial z_{n}}
\right)$ is the gradient of $f$, and the last condition just expresses that 
each entry of the (row) vector $(\nabla f)A$ is divisible by $f$ in $P$. 
The columns of $A$ can then be viewed as the coefficients of a basis, with 
respect to the partial 
derivatives $\partial/\partial z_{i}$,  of the 
logarithmic vectorfields along the divisor $f=0$.\eop
\end{proposition}

\ni The normal crossing divisor 
$D= \{x_1\cd x_n=0\}$ provides a simple example: Saito's criterion shows 
that the vector fields $x_1\p/\p x_1,
\ld,x_n\p/\p x_n$ form a 
 basis for $\Der(-\log D)$. This free
divisor has the striking property that $\Der(-\log D)$ has a 
basis consisting of vector fields 
that are homogeneous 
 of weight zero
with respect to the natural grading.
 Among free hyperplane
arrangements it is the only one 
with this property
(\cite{orlikterao} Chapter 4). Until recently,
the only other free divisor with this property known to either of the
authors of this paper was the ``{\em bracelet\/}'', 
the discriminant in the space of 
binary cubics (see
\cite{cg}, and \cite{mwik}, 
where it is described in some detail, though
not under this name). 
\begin{definition} The free divisor $D$ is {\em linear} if $\Der(-\log D)$
has a 
basis consisting of vector fields of weight zero --- that is, 
all of whose coefficients are linear functions of the variables.
\end{definition} 
\ni Here we show that far from being uncommon, linear free divisor are abundant.
We show that the set of degenerate, or non-generic, orbits in the
representation space $\Rep(Q,\d)$ of a quiver with dimension vector $\d$,
is a linear free divisor
whenever $\d$ is a real Schur root (definition in Section \ref{background}) 
of $Q$, and provided that a natural
condition on the existence of ``codimension 1'' degeneracies holds -
a condition which is closely related to Damon's condition on the
existence of ``Morse-type singularities'' mentioned above.\\

\ni Since we hope that our paper will be read by singularity theorists, we 
include some background on quiver representations.
\section{Linear free divisors}
Suppose that $D$ is a linear free divisor, and let $\chi_1,\ld,\chi_n$
be a 
basis consisting of weight-zero vector fields. Since 
the weight of the Lie bracket of any two homogeneous vector fields is
the sum of their weights, $\chi_1,\ld,\chi_n$ 
form the basis of
an $n$-dimensional
Lie algebra $L_D$ over $k$, as well as 
a basis of
the free $\OO$-module $\Der(-\log D)$.
Consider the standard action of $\Gl_n(k)$ on $k^n$. The vector field
$x_i\p/\p x_j$ is the infinitesimal generator of this action corresponding to
the elementary matrix $E_{ij}$ (1 in the $i$-th row and $j$-th column,
zeroes elsewhere). It follows that $L_D$ is the image, under the infinitesimal
action, of an $n$-dimensional Lie subalgebra of $\gl_n(k)$, which we denote
$\gg_D$. In the complex case, if the exponential of $\gg_D$ 
is a closed subgroup $G_D$
of $\Gl_n(\CC)$, then $G_D$ has an open orbit in $\CC^n$ 
and $D$ is its complement. 
This follows easily from Mather's lemma on
Lie group actions
(\cite{matherIV} Lemma 3.1), which gives sufficient conditions for a connected
submanifold of a manifold to lie in a single orbit of the action of a Lie group
$G$: that\\

(i) at each point of $X$, $T_xX$ should be contained in the tangent space
to the $G$ orbit of $x$, and 

(ii) the dimension of this orbit 
should be constant for $x\in X$.\\

\ni Taking $X=\CC^n\setminus D$, both conditions evidently hold here.\\

\ni In all examples known, this indeed applies.
To find linear free divisors one 
may thus look
for $n$--dimensional
Lie groups acting on $k^n$ with an open orbit. It is precisely these that
the representation theory of quivers offers in abundance.
Indeed, in that situation, the Lie groups $G_{D}$ are reductive. Examples
of some nonreductive groups that also give rise to linear free divisors 
will be presented in \cite{gmns}. Here we mention just one series.
\begin{example}{\em  The group $B_n(k)$ of upper triangular $n\times n$ matrices
acts on the space $\Sym_n(k)$ of symmetric matrices by 
$$B\cdot S = {}^tB\,S\,B.$$
There is an open orbit; the equation of the complement is the product of
$n$ nested symmetric determinants, beginning with the top left hand entry
($1\times 1$ determinant) in the symmetric matrix $S$ and continuing
with the determinant of the top left hand $2\times 2$ block, the 
determinant of the top left hand $3\times 3$ block, etc.
}
\end{example}

\section{Representations of Quivers}\label{background}
A {\it quiver } is a 
finite
directed graph. That is, it consists of a 
finite
set
$Q_0$ of nodes (or vertices), 
and a finite set of arrows $Q_{1}$ equipped with two maps $h,t:Q_{1}\to Q_{0}$ 
that assign to each arrow $\pp\in Q_1$ its head $h\pp$ and tail $t\pp$ in 
$Q_0$.
A {\em representation} $V$ of a quiver 
$Q$ consists of a choice of vector space $V_x$
for each node $x$, and a $k$-linear map $V(\pp):V_{t\pp}\to V_{h\pp}$ 
for each arrow $\pp\in Q_1$. 
The representation is {\em finite dimensional\/} if each 
$V_{x}$ is a finite dimensional vector space. 

If $W$ is a second such representation, then a {\em morphism of
representations\/} $\psi:W\to V$ is a family of $k$--linear maps 
$\psi_{x}:W_{x}\to V_{x}, x\in Q_{0},$ such that for each 
$\pp\in Q_{1}$ the square
$$
\xymatrix{
W_{t\pp}\ar[rr]^{W(\pp)}
\ar[d]_{\psi_{t\pp}}&&W_{h\pp}
\ar[d]^{\psi_{h\pp}}\\
V_{t\pp}\ar[rr]^{V(\pp)}&&V_{h\pp}
}
$$
commutes. 
The $k$--vector space of all morphisms of representations from 
$W$ to $V$ is denoted $\Hom_{Q}(W,V)$. The so-defined category of 
(finite dimensional) representations of $Q$ is {\em abelian\/}. 
Moreover, it is {\em hereditary\/}, which means that the extension 
groups in this abelian category --- denoted $\Ext^{i}_{Q}(W,V)$, or 
$\Ext^{i}_{kQ}(W,V)$ if we wish to specify the coefficients --- 
vanish whenever $i\ge 2$.

Once we fix the dimensions of the spaces at
each node, by assigning to $Q$ a {\it dimension vector} $\d\in\NN^{Q_0}$, 
we can consider the $k$-vector space of representations 
$$\Rep(Q,\d)=\prod_{\pp\in Q_1}\Hom_k(V_{t\pp},V_{h\pp})
\simeq \prod_{\pp\in Q_1}\Hom_k\bigl(k^{\d(t\pp)},k^{\d(h\pp)}\bigr).$$
The group $\Gl(Q,\d)=\prod_{x\in Q_0}\Gl_{\d(x)}(k)$ acts on 
$\Rep(Q,\d)$ by 
$$(g_x)_{x\in Q_0}\cdot (V(\pp))_{\pp\in Q_1}=\bigl(
g_{h\pp}\circ V(\pp)\circ
g_{t\pp}^{-1}\bigr)_{\pp\in Q_1}.$$
The orbits of this group action are the isomorphism classes of 
$Q$--representations with the prescribed dimension vector.
It will be from this action that we obtain the generators of $\Der(-\log D)$
for the linear free divisors we construct.
\\

Given $V'\in\Rep(Q,\d')$ and 
$V''\in\Rep(Q,\d'')$, the {\it direct sum} $V'\oplus V''
\in\Rep(Q,\d'+\d'')$ is the representation with 
$(V'\oplus V'')_x=V'_x\oplus V''_x$ for $x\in Q_0$ and
$$(V'\oplus V'')(\pp)=
\left(\begin{array}{cc} V'(\pp)&0\\0&V''(\pp)\end{array}
\right).$$
A given representation $V\in\Rep(Q,\d)$ is {\it decomposable} if it is
the direct sum of subrepresentations --- that is, if there are representations
$V'\in\Rep(Q,\d')$ and $V''\in\Rep(Q,\d'')$ such that $V=V'\oplus V''$.
In this case, of course, $\d=\d'+\d''$.

A quiver $Q$ is a {\it Dynkin quiver} if the underlying undirected graph
$\overline Q$
is  
a disjoint union of Dynkin diagrams
of type $A_n$, $D_n$, $E_6, E_7$ or $E_8$. Dynkin quivers
are ubiquitous in the theory of representations of quivers, and central in
this paper.
\begin{example}{\em Let $Q$ be the 
Dynkin quiver of type $A_3$
$$
\xymatrix{\bullet\ar[r]^{\displaystyle A}&\bullet\ar[r]^{\displaystyle B}&\bullet}
$$

\ni (i) With dimension vector $(1,1,1)$ any representation in which
each of the morphisms is non-zero is indecomposable. \\
(ii) Indeed, these are the only
indecomposable representations 
whose dimension vector is {\em sincere\/}, 
meaning that it is nonzero at each node. For example,
if $\d=(1,2,1)$
there is no indecomposable representation. 
Representations in which $BA\neq 0$ decompose as the direct sum
$$k\st{A}{\too}\im\,A \cong k \xrightarrow{B|_{\text{im}\,A}} 
k\quad\oplus\quad 0\too\ker\,B\cong k\too 0$$
Representations in which $BA=0$ and $A\neq 0$ decompose as
$$k\st{A}{\too}\im\,A\cong k \too 0\quad\oplus\quad 0\too k^{2}/\im\,
A\cong k\st{B}{\too}k$$
where the middle term in the second summand 
can be viewed as a complement to $\im\,A$.
Representations in which $A=0$ decompose as
$$k\too 0 \too 0\quad\oplus\quad 0\too k^2\st{B}{\too}k.$$

Similarly,  any representation with $\d=(l,m,n)$ for $l,m,n\ge 0$ decomposes as:
\begin{alignat*}{4}
(k\too &0\too 0)^{\oplus a}&\quad\oplus\quad&(0\too k&\too 0)^{\oplus b}&\quad\oplus\quad&(0\too &0\too k)^{\oplus c}\\
\oplus\quad (k\st{1}\too &k\too 0)^{\oplus d}&\quad\oplus\quad&(0\too k&\st{1}\too k)^{\oplus e}&\quad\oplus\quad&(k\st{1}\too &k\st{1}\too k)^{\oplus f}\,,
\end{alignat*}
where 
\begin{align*}
a &=\dim\,\ker A\quad,\quad b = \dim\,\ker B/(\im\,A\cap\ker\,B)\quad,\quad c=\dim\,\ \text{cok}\,B\,,\\
d &=\dim\, \ker BA/\ker A\quad,\quad
e = \dim\,\ \im\,B/\im\,BA\,,\\
 f &= \dim\,\ \im\,BA = l-a-d = m - b - d- e = n - c - e\,. 
\end{align*}

}
\end{example}

\begin{definition}\label{root}
The dimension vector $\d$ is a {\em root} of $Q$ if 
$\Rep(Q,\d)$ contains an
indecomposable representation. The root is {\em real\/} if $\Rep(Q,\d)$ 
contains exactly one orbit of, necessarily isomorphic, indecomposable 
representations. It is {\em imaginary\/} if there is a family of non-isomorphic 
indecomposable representations. If a general representation in $\Rep(Q,\d)$ 
is indecomposable, then $\d$ is a {\em Schur\/} root.\footnote{Some authors 
define a Schur root as a root $\d$
for which $\Rep(Q,\d)$ contains a `brick' --- a representation $V$ for which
$\End_Q(V)=k$. If $\d$ is a Schur root in this sense, then by the upper 
semicontinuity of $\dim\,\End_Q(V)$ with respect to $V$, the general 
representation also has 
endomorphism ring $k$, and so is indecomposable. Conversely, 2.7 of \cite{krre}
shows that if the general representation is indecomposable then it is a brick.
So the two versions of the definition are equivalent.} 

{\em The frequent use of the term ``root'' in these definitions is no 
coincidence,
as we will see below.}
\end{definition}

A crucial role in the representation theory of quivers is played by
the {\it Euler form}, a bilinear form on the space $\NN^{Q_0}$ of 
dimension vectors. It is defined by
$$
\langle \mathbf{e}, \mathbf{d}\rangle 
=\sum_{x\in Q_0}e_xd_x-\sum_{\pp\in Q_1}e_{t\pp}d_{h\pp}
=
\dim\,\prod_{x\in Q_0}\Hom(W_x,V_x)-
\dim\,\prod_{\pp\in Q_1}\Hom_k(W_{t\pp},V_{h\pp})
$$
for any $W\in R(Q,\e)$ and $V\in R(Q,\d)$, and accordingly we sometimes denote $\la \e,\d\ra$ by
$\la W, V\ra$.

The {\it Tits form} on the space of dimension vectors is the associated
quadratic form, $q(\d)=\la\,\d,\d\ra$. Observe that 
the Tits form does not depend on the orientation of the arrows. Indeed, 
it is used to calculate the members of the root system of the Kac--Moody 
Lie algebra attached to the underlying graph $\overline Q$, and those 
roots with nonnegative components are precisely the roots for $Q$, 
regardless of the orientation of the arrows, see \cite{kac}. 
For example, if $\overline Q$ is a Dynkin diagram 
then $\d\in \QQ^{|Q_0|}$ is a root of the corresponding semi-simple 
Lie algebra, in the classical sense, if and only if
$q(\d)=1$. In particular, all roots are real in this case.

Choosing an ordering of the nodes in $Q_{0}$, 
we may write $\langle\e,\d\rangle =  
\e E\d^{T}$, where $\mathbf{e}, \mathbf{d}$ 
are thought of as row vectors, and $E$ is the corresponding 
{\em Euler matrix\/}. Its entries are 
$E_{x,y} = \delta^{y}_{x}-
\#\{\pp\in Q_{1}\mid t\pp = x, h\pp = y\}$, with 
$\delta^{y}_{x}$ denoting the Kronecker delta. Put differently, 
$E= I_{|Q_{0}|} - A$, where
 $I_{|Q_{0}|}$ is the identity matrix of 
the indicated size and the matrix entry $A_{x,y}$ records the 
number of arrows from $x$ to $y$ in $Q_{1}$. The matrix associated 
to the Tits form is then $C=E+E^{T}$, the {\em Cartan matrix\/} of $Q$, 
which coincides with the usual Cartan matrix of the associated Dynkin 
diagram $\overline Q$, in case $Q$ is a Dynkin quiver\footnote{More generally, 
$C$ is the Cartan matrix of the Kac--Moody Lie algebra associated to 
$\overline Q$, for an arbitrary finite quiver $Q$ without oriented cycles, 
see \cite{kac} again.}.

The following simple result is useful for the actual calculation of 
the linear free divisors below.
\begin{lemma}
\label{lem:einv}
If $Q$ is a finite quiver without oriented cycles, 
then its Euler matrix is {\em invertible\/}. The inverse is given by
$
E^{-1} = I_{|Q_{0}|} + A'
$,
where $A'_{x,y}$ equals the number of {\em directed paths\/} from $x$ to $y$.
%
\eop
\end{lemma}
\ni
Now we recall the trichotomy of the representation 
theory of quivers:

\begin{definition}
\ni A quiver $Q$ is of {\em finite representation type} if $Q$ has only finitely
many indecomposable representations, up to isomorphism. The quiver is 
{\em wild} if its representation theory is at least as complicated as that
of the quiver
$$
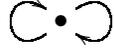
\xymatrix@R=0.8in@C=0.8in{\bullet\ar@(ur,dr)\ar@(dl,ul)}
$$
\vskip 10pt
\ni The quiver is 
{\em tame} if it is neither of finite representation type, nor wild\,\footnote{
The reader should be aware that the definition often is ``tame''=``not wild'', 
thus, different from our usage here. 
}. 
\end{definition}

Gabriel (\cite{gabriel},\cite{gabriel1}) showed
\begin{theorem}
\label{thm:gabriel} 
A connected quiver $Q$ is of finite representation
type if and only if it is a Dynkin quiver. 
Assigning to an isomorphism class of indecomposable representations of
$Q$ its dimension vector induces then a
bijection between these classes and the 
positive roots of the underlying Dynkin diagram.
\eop
\end{theorem}
\ni
The last part of this result can be restated thus: 
if $\d$ is a positive root of the underlying Dynkin diagram
(as listed, for example, in the appendix to \cite{B}) then $\d$ is also
a root of any associated Dynkin quiver $Q$, in the sense of Definition 
\ref{root}.
Moreover, in this case each 
root is a real Schur root: there is a (unique) open orbit
in $\Rep(Q,\d)$ whose points correspond to indecomposable representations. 
A good account of all this can be found in \cite{bgp}.\\
The class of tame quivers has a similar characterisation: 
\begin{theorem}{\em (\cite{df},\cite{naz})} 
A connected quiver is tame if and only if the underlying
undirected graph is an extended Dynkin diagram. \eop
\end{theorem}
\ni
Finally, in what follows we will need a result of V.Kac (\cite{kac}):
\begin{proposition}\label{kacs} 
Let $Q$ be a connected quiver whose proper subquivers are
all either of finite or tame type. Then 
a dimension vector 
$\d$ is a real 
root if and only if
$q(\d)=1$, and it is an imaginary 
root if and only if $q(\d)\leq 0$.\eop
\end{proposition}
\section{The fundamental exact sequence}
Let $V$ and $W$ be representations of the quiver $Q$. In \cite{ringelalg},
Ringel introduced the following exact sequence $\s E^W_V$ of vector spaces:
\begin{align}
\label{fes}
&0\to \Hom_Q(W,V)\to \prod_{x\in Q_0}\Hom_{k}(W_x,V_x)\st{d^W_V}{\too}
\prod_{\pp\in Q_1}\Hom_k\bigl(W(t(\pp)),V(h(\pp))\bigr)\st{e^W_V}{\too} 
\Ext^1_{Q}(W,V)\to 0\,.
\end{align}
The morphism $d^W_V$ is defined by 
$$d^W_V\bigl((\psi_{x})_{x\in Q_0}\bigr)=
\bigl(\psi_{h(\pp)}\cc W(\pp)
-V(\pp)\circ \psi_{t(\pp)}\bigr)_{\pp\in Q_1};$$
the component of $d_V^W((\psi_{x}))$ corresponding to $\pp\in Q_1$
measures non-commutativity of the diagram
$$
\xymatrix{
W_{t\pp}\ar[rr]^{W(\pp)}
\ar[d]_{\psi_{t\pp}}&&W_{h\pp}
\ar[d]^{\psi_{h\pp}}\\
V_{t\pp}\ar[rr]^{V(\pp)}&&V_{h\pp}
}
$$
whence it is clear that $\ker d^W_V$ is indeed equal to $\Hom_Q(W,V)$.\\

\ni To define $e^W_V$, 
from $\theta = (\theta_\pp)_{\pp\in Q_1}$ we construct a new representation
$Z$ of $Q$ and an exact sequence, 
\begin{align*}
e^W_V(\theta)\quad\equiv\quad 0\to V\st{i}{\too} Z\st{j}{\too} W\to 0\,,
\end{align*}
by the following recipe:
$Z_x=V_x\oplus W_x$ for each $x\in Q_0$, $i_x:V_x\to V_x\oplus W_x$
and $j_x:V_x\oplus W_x\to W_x$ are the standard inclusion and projection,
and for each $\pp\in Q_1$,
$Z(\pp):V_{t\pp}\oplus W_{t\pp}\to V_{h\pp}\oplus
W_{h\pp}$ has matrix
$$\left(
\begin{array}{cc}
W(\pp)&\theta_{\pp}\\
0&V(\pp)
\end{array}
\right).
$$ 
It is straightforward to check that the short exact sequence $e^W_V(\theta)$
of representations
of $Q$ is split if and only if $\theta=d^W_V(\psi)$ for some $\psi
\in\prod_{x\in Q_0}\Hom(W_x,V_x)$, and that $e^W_V$ is onto.

Exactness of 
the sequence $\s E^W_V$
implies that
$$\la \e,\d\ra =\dim_k\Hom_{Q}(W,V) - \dim_k\Ext^1_{Q}(W,V)$$
 for any $W\in R(Q,\e)$ and $V\in R(Q,\d)$,
 so that the expression on the right hand side depends only on the dimension
vectors and not on the choice of representations, although evidently the
dimensions of 
$\Ext^1_{Q}(W,V)$ and $\Hom_Q(W,V)$ do depend on the choice of $V\in
R(Q,\d)$ and $W\in R(Q,\e)$.\\

\ni The fundamental sequence $\s E^W_V$ plays two roles in what follows. 
In the next section we show how 
to reinterpret it in terms of the deformation theory of
representations, where it may become more familiar to singularity-theorists.
From this we will see how free divisors appear naturally in this context.
\\

\ni Second, following Schofield \cite{schofield}, 
we use it
to generate semi-invariants of the representation space $R(Q,\mathbf{d})
=\prod_{\pp\in Q_1}\Hom_k(k^{d(t\pp)},k^{d(h\pp)})$, 
and thereby find explicit
equations for the free divisors, in Sections \ref{eqns} and \ref{e8}. 
\section{Deformations of representations}\label{defs}
Recall that the group $\Gl(Q,\d)$ acts on $\Rep(Q,\d)$ by
$$(g_x)_{x\in Q_0}\cdot (V(\pp))_{\pp\in Q_1}=\bigl(
g_{h\pp}\circ V(\pp)\circ
g_{t\pp}^{-1}\bigr)_{\pp\in Q_1}.$$
The orbit of $V$ in $\Rep(Q,\d)$ is open if and only if the associated map
$$\a_V:\Gl(Q,\d)\to \Rep(Q,\d)$$
sending $g$ to $g\cdot V$ is a submersion, and for this it is enough
that it be a submersion at the identity. The Lie algebra $\gl(Q,\d)$ of
$\Gl(Q,\d)$ is
$$\prod_{x\in Q_0}\End(k^{\d(x)})=\prod_{x\in Q_0}\Hom(k^{\d(x)},k^{\d(x)}),$$
and the tangent space to $\Rep(Q,\d)$ at $V$ is $\Rep(Q,\d)$ itself, i.e.
$\prod_{x\in Q_0}\Hom_k(k^{\d(t\pp)},k^{\d(h\pp)})$. The derivative
of $\a_V$ 
at the identity in $\Gl(Q,\d)$ is
precisely the map $d^V_V$ of the exact sequence
 $\s E^V_V$.
In fact we may canonically identify $\Ext^1_{Q}(V,V)$ with 
$T^1(V)$ for the associated deformation theory, though we will not
make any formal use of this identification.
\\

\ni A deformation, in the analytic category, of a representation $V$ is, 
by definition, the germ of 
an analytic map $(B,0)\to (\Rep(Q,\d),V)$. 
If $(B,0)$ is smooth,
a deformation $\s V:(B,0)\to
(\Rep(Q,\d),V)$ is {\em versal\/} 
if and only if it is {\em complete\/}, that is, 
if every
other deformation $\s V':(B',0)\to (\Rep(Q,\d),V)$ is 
equivalent to one induced 
from it by base-change $\eta:(B',0)\to (B,0)$. 
The equivalence here is the existence
of a map-germ $g:(B,0)\to (\Gl(Q,\d),1)$ such that 
$$\s V'(b')=g(b')\cdot \s V(\eta(b')).$$
Thus it is evident that $\Rep(Q,\d)$
itself, or more precisely the identity map $(\Rep(Q,\d),V)\to (\Rep(Q,\d),V)$,
is a versal deformation; for any other deformation $\s V'$,
the base change map $\eta$ is simply $\s V'$
itself, and $g$ is the constant map taking the value $1$. 
The slice theorem from the theory of smooth group
actions is now enough to establish the versality of any deformation obtained
from this one by restricting its domain to any smooth space-germ transverse
to the orbit of $V$, or indeed by pulling it back by any map-germ $(B,0)\to
(\Rep(Q,\d),V)$ transverse to the
orbit of $V$.
These considerations imply the Artin--Voigts Lemma: 
that the dimension of $\Ext^1_{Q}(V,V)\cong T^{1}(V)$ equals the 
codimension of the orbit of $V$ in $\Rep(Q,\d)$. In particular, if 
there is an open orbit, then the representations therein have no 
self-extensions: they are {\em rigid\/} as representations.
\\

\ni Now we consider the relative $T^1$, obtained by regarding the coefficients
of the morphisms $V(\pp)$ as variables. This can be done in the analytic,
formal or algebraic category, and amounts to no more than tensoring the
exact sequence 
 $\s E^V_V$
with the appropriate ring, or sheaf, of functions 
--- $\OO_{\Rep(Q,\d)}$, 
$k[\Rep(Q,\d)^{*}]$ or $k[[\Rep(Q,\d)^{*}]]$. 
We refer to these indistinctly as $R$. 
The module (sheaf) of vector fields on $\Rep(Q,\d)$ is 
$\theta_{R} = \Der_k(R)\cong \Rep(Q,\d)\otimes_{k}R$, 
and the $k$-linear
map $\mathfrak{gl}(Q,\d))\to \Rep(Q,\d)$ extends to a morphism of $R$-modules
$\gl(Q,\d)\otimes_{k}R\to\theta_R$ whose cokernel can be viewed both as
$\Ext^1_{RQ}(M,M)$ for the universal representation $M$ 
of the quiver $Q$ with coefficients in $R$, 
and as the relative $T^1$ of the
versal deformation $\mathfrak{i}:\Rep(Q,\d)\to \Rep(Q,\d)$, 
denoted 
$T^1(\mathfrak{i}/\Rep(Q,\d))$. The surjection 
$\theta_R\to T^1(\mathfrak{i}/\Rep(Q,\d))$ 
is the {\em Kodaira-Spencer map\/} of the versal
deformation $\mathfrak{i}$.

The kernel of this projection is the space of simultaneous endomorphisms of
the representations $V\in\Rep(Q,\d)$, or, in other words, the endomorphism
ring of the universal representation $M$. Provided the general representation
in $\Rep(Q,\d)$ is indecomposable, this ring is isomorphic to 
$R$. 
Let us understand why this is so. It is 
clear that if $V\in \Rep(Q,\d)$ is any representation then for 
each $\l\in k^*$ we have $(\l I_{d_x})_{x\in Q_0}\in \Aut_k(V)$, and similarly
$(\l I_{d_x})_{x\in Q_0}\in \End_Q(V)$ for $\l\in k$. 
If $V$ is stably indecomposable (that is, if there
is a neighbourhood of $V$ in $\Rep(Q,\d)$ consisting of indecomposable 
representations) 
then this copy of $k$ accounts for all of $\End_Q(V)$ 
(see e.g. \cite{krre}, 2.7 ).
Now if the general representation in $\Rep(Q,\d)$ is indecomposable ---
which means that
$\d$ is a Schur root ---
then  at each of these representations, any endomorphism
of the universal representation $M$ must be a scalar. By density, the
same must be true everywhere, and so $\End_R(M)$ can be identified with 
$R$.

The cokernel of 
the
inclusion of Lie algebras $0\to k\to \mathfrak{gl}(Q,\d)$
is, by definition, $\mathfrak{pgl}(Q,\d)$, and we can identify the cokernel of 
the inclusion of free $R$-modules 
$0\to R\to \mathfrak{gl}(Q,\d))\otimes R$ with 
$\mathfrak{pgl}(Q,\d)\otimes R$. Thus, provided the generic representation
in $\Rep(Q,\d)$ is indecomposable, we have a short exact sequence
\beq\label{fes2}
0\to\mathfrak{pgl}(Q,\d)\otimes_kR\st{\tilde d^M_M}{\too} 
\theta_R\to \Ext^1_{RQ}(M,M)\to 0.
\eeq
Even without generic indecomposability, we still have an exact sequence
\beq\label{fes3}
\mathfrak{pgl}(Q,\d)\otimes_kR\st{\tilde d^M_M}{\too} 
\theta_R\to \Ext^1_{RQ}(M,M)\to 0.
\eeq
Let $D$ be the support of $\Ext^1_{RQ}(M,M)=T^1(\mathfrak{i}/\Rep(Q,\d))$,
with (possibly non-reduced) coordinate ring 
$R[D]=R/\s F_0\bigl(\Ext^1_{RQ}(M,M)\bigr),$ where $\s F_0$ means
zero'th Fitting ideal. 

\begin{proposition}\label{fp} 

\ni {\em (i)\/} $D$ is the set of non-rigid 
representations. 
Its open complement is the set of rigid
representations\footnote{
Singularity theorists might prefer 
`stable' to `rigid'; 
however in representation theory the term `stable' often
refers to its meaning in geometric invariant theory, so here we use `rigid'.
}.
\\
\ni If $q(\d)=1$ and the general representation in
{\em $\Rep(Q,\d)$\/} is indecomposable, 
thus, $\d$ is a Schur root,
then\\

\ni  {\em (ii)\/} $D$ is a divisor in $\Rep(Q,\d)$.\\

\ni {\em (iii)\/} {\em $\Ext^1_{RQ}(M,M)$} is a maximal Cohen-Macaulay 
$R[D]$-module.\\ 

\ni {\em (iv)\/} The image of 
$\tilde d^M_M:\mathfrak{pgl}(Q,\d)\otimes_kR\to \theta_R$ is 
contained in $\Der(-\log D)$.\\

\end{proposition} 
\ni\Proof (i) Let $m_V$ be the maximal ideal of $R$ corresponding to 
$V\in\Rep(Q,\d)$. By right-exactness of tensor product, tensoring 
the sequence
(\ref{fes3}) with $R/m_V$ gives the exact sequence
$$\mathfrak{pgl}(Q,\d)\st{ d^V_V}{\too} 
\Rep(Q,\d)\to \Ext^1_Q(V,V)=T^1(V)\to 0.
$$
This establishes (i).\\

\ni (ii) Since now $\mathfrak{pgl}(Q,\d)\otimes_kR$ and $\theta_R$ are free
$R$-modules of the same rank, $\s F_0(\Ext^1_{RQ}(M,M))$ is generated by
$\det(\tilde d^M_M)$, and so $D=\supp(\Ext^1_{RQ}(M,M))=
V\bigl(\s F_0(\Ext^1_{RQ}(M,M))\bigr)
=V(\det(\tilde d^M_M))$.
\\

\ni (iii) Exactness of the sequence (\ref{fes2}) implies, by the 
Auslander-Buchsbaum formula,
that 
$$\depth_R(\Ext^1_{RQ}(M,M)) = \dim\,R-1 = \dim\,\Ext^1_{RQ}(M,M)\,,
$$
where ``$\dim$'' here refers to Krull dimension.
 
Hence 
$\Ext^1_{RQ}(M,M)$ is a Cohen-Macaulay $R$-module. It is annihilated by 
$\s F_0\bigl(\Ext^1_{RQ}(M,M)\bigr)$, so is an $R[D]$-module, and as such,
a maximal Cohen-Macaulay 
module.

\ni (iv) The vector fields in 
$\tilde d^M_M\bigl(\pgl(Q,\d)\otimes_k R\bigr)$ 
are 
infinitesimal generators of the action of $\Gl(Q,\d)$ on $\Rep(Q,\d)$,
and are thus tangent to all its orbits. So they are tangent to
$D$, which is a union of orbits.\\
\eop

\ni Note that by \ref{thm:gabriel}, if $Q$ is 
a Dynkin quiver then in order for
(ii)-(iv) to hold we need only require that $q(\d)=1$.\\

\ni If the conditions of \ref{fp}(ii)-(iv) hold, and moreover the
vector fields in 
$\tilde d^M_M\bigl(\pgl(Q,\d)\otimes_kR\bigr)$
{\it generate} $\Der(-\log D)$ then $D$ is a linear free divisor, since 
by exactness of (\ref{fes2}), $\Der(-\log D)$ is free over $R$.
Saito's criterion (\ref{scrit} above) shows that in order that they 
do generate, it is enough that 
$\det\bigl(\tilde d^M_M\bigr)$ be reduced. \\

\ni Thus, we obtain the following result: 
\begin{corollary} With the conditions of {\em \ref{fp}(ii)-(iv)}, suppose in 
addition that $D$ is reduced. Then it is a linear free divisor.\eop
\end{corollary}

\ni From now on we will refer to the divisor $D$ of non-rigid representations
in $\Rep(Q,\d)$ as the {\it discriminant}
and call $\Delta = \det\bigl(\tilde d^M_M\bigr)$ its {\em canonical equation}.\\

\ni Suppose that $D$ is reduced at $V$. Then by Saito's criterion, the 
vector fields in 
$\tilde d^M_M\bigl(\pgl(Q,\d)\otimes_kR\bigr)$ generate
the stalk at $V$ of the sheaf $\Der(-\log D)$. If $V$ is a regular point
of $D$, then the tangent space $T_VD$ is equal
to $d^V_V\bigl(\pgl(Q,\d)\bigr)\subseteq \Rep(Q,\d)$. 
It follows that the deformation
of $V$ obtained by following any smooth curve transverse to $D$ is
versal. For the same reason, any deformation in a direction
tangent to $D$ is infinitesimally trivial at $V$. Since the same holds at
any nearby point, any deformation of $V$ in the smooth part of $D$ is globally
trivial.
In terms of the group action, this means that then each 
irreducible component of $D$ contains 
a dense
open orbit, and for each 
representation $V$ in such an orbit, $T^{1}(V)$ will be one-dimensional. 
We now investigate further the relation between the dimension of $T^{1}(V)$ 
for a generic representation on such a component and the multiplicity with 
which  that component occurs in the discriminant.

\begin{lemma} Let $D_{j}$ be an irreducible component of $D$, and $h_{j}$ its 
reduced equation, $m_j$ the multiplicity of $h_{j}$ in 
$\det \bigl(\tilde d^M_M\bigr)$, and $V_{j}$ a generic representation on 
$D_{j}$. One has then
$m_{j} \ge \dim_kT^1(V_j)$ and equality holds if and only if 
$h_{j}$ annihilates $\Ext^{1}_{RQ}(M,M)$. In particular, 
$m_{j}=1$ forces $\dim\, T^{1}(V_{j})$ to be one-dimensional and the 
orbit generated by $V_{j}$ to be dense in $D_{j}$. 
\end{lemma}
\ni\Proof 
Let $p$ be the ideal $(h_j)$. Then the localisation $R_{p}$ is a discrete 
valuation ring. We must have 
\beq\label{meal}
\Ext^1_{RQ}(M,M)\otimes_RR_{p}\simeq\bigoplus_1^\ell R_p/(pR_p)^{\alpha_{t}}
\eeq 
for some positive integers $\alpha_{t}$; it follows that the matrix 
$\tilde d^M_M$
is equivalent, over $R_p$, to a matrix of the form
$\text{diag}(h_j^{\alpha_{1}},..., h_{j}^{\alpha_{\ell}})\oplus I_{r-\ell}$, 
a block matrix formed of the indicated diagonal matrix and the identity 
matrix $I_{r-\ell}$, where $r=\dim_k\pgl(Q,\d)$. Evidently $\det \tilde 
d^M_M=(h_j)^{\sum_{t=1}^{\ell}\alpha_{t}}$
in $R_p$, and so $\sum_{t=1}^{\ell}\alpha_{t}=m_j$. Moreover, by 
(\ref{meal}), $\sum_{t=1}^{\ell}\alpha_{t}$ is also
the rank of $\Ext^1_{RQ}(M,M)$ at a generic point $V_j$ of $D_j$. Dividing by
the maximal ideal 
$m_{V_j}$, we see that then $\ell$ is equal to $\dim_k\Ext^1_Q(V_{j},V_{j})$.
Therefore, $m_{j}= \sum_{t=1}^{\ell}\alpha_{t} \ge \ell = \dim_k\Ext^1_Q(V_{j},
V_{j}) = \dim T^{1}(V_{j})$. Clearly, $m_{j}=\ell$ if and only if each 
$\alpha_{t}=1$ if and only if $h_{j}$ annihilates $\Ext^1_{RQ}(M,M)$. \eop
In the case of Dynkin quivers, it follows that the discriminant is 
indeed reduced, as we show next.

\begin{proposition}
Let $\d$ be a real Schur root of a Dynkin quiver $Q$ and assume that 
$V\in \Rep(Q,\d)$ satisfies $\dim\, T^{1}(V) =1$.
If $D'\subseteq D$ denotes the irreducible component of the discriminant that 
contains $V$ and $h'=0$ is its reduced equation, then $h'$ divides 
$\Delta = \det \bigl(\tilde d^M_M\bigr)$ with multiplicity one.
\end{proposition}

\ni\Proof We begin by clarifying in general what it means that $D'$ 
appears with multiplicity one, if we know already that the generic 
representation on it has one-dimensional $T^{1}$: As $T^{1}(V)=\Ext^{1}_{Q}(V,V)$ 
is one-dimensional, the semi-universal deformation of $V$ as a representation 
of $Q$ has a one-dimensional base. Because $V$ deforms into a rigid 
representation generically, its reduced discriminant consists just of the origin. 
By Openess of Versality, it suffices to prove that the discriminant in that 
semi-universal deformation is indeed reduced. If $\mathfrak V$ is the universal 
module over $k[\![t]\!]$, the (formal) base ring of the semi-universal 
deformation, it suffices to show that $\Ext^{1}_{k[\![t]\!]Q}({\mathfrak V},
{\mathfrak V})$ is a one-dimensional vector space. Now $\Ext^{1}_{k[\![t]\!]Q}
({\mathfrak V},{\mathfrak V})$ is concentrated on the discriminant, thus a 
finite dimensional vector space. Moreover, $\Ext^{1}_{k[\![t]\!]Q}({\mathfrak V},
{\mathfrak V})\otimes_{k[\![t]\!]}k \cong
\Ext^{1}_{Q}(V,V)\cong k$, whence as $k[\![t]\!]$--module 
$\Ext^{1}_{k[\![t]\!]Q}({\mathfrak V},{\mathfrak V})\cong k[\![t]\!]/(t^{m})$ 
for some $m$. We need to show that $m=1$, and this can be achieved by 
establishing that the following natural projection, in its various guises:
\begin{align*}
\xymatrix{
\Ext^{1}_{k[\![t]\!]Q}({\mathfrak V},{\mathfrak V})\otimes_{k[\![t]\!]}
k[\![t]\!]/(t^{2})\ar[d]& \cong& 
\Ext^{1}_{k[\![t]\!]Q}({\mathfrak V},{\mathfrak V}/(t^{2}){\mathfrak V})
\ar[d]& \cong&
k[\![t]\!]/(t^{m}, t^{2})\ar[d]\\
\Ext^{1}_{k[\![t]\!]Q}({\mathfrak V},{\mathfrak V})\otimes_{k[\![t]\!]}k 
&\cong& \Ext^{1}_{k[\![t]\!]Q}({\mathfrak V},{\mathfrak V}/(t){\mathfrak V}) 
\cong \Ext^{1}_{kQ}({V},{V})&\cong& k
}
\end{align*}
is an isomorphism.
 To this end, let
\begin{align}
\label{al:1}
0\to V \xrightarrow{i} W \xrightarrow{p} V\to 0
\end{align}
represent a nontrivial element in the one-dimensional vector space 
$\Ext^{1}_{Q}(V,V)$. Define an action of $t$ on $W$ through $t(w) = ip(w)$. 
Clearly, $t^{2} = ipip =0$ on $W$, whence the $Q$--representation $W$ 
becomes as well a $k[\![t]\!]/(t^{2})$--module. Infinitesimal deformation 
theory says that indeed $W\cong {\mathfrak V}/t^{2}{\mathfrak V}$, and 
that the extension above can be viewed as an extension of $k[\![t]\!]$--modules,
\begin{align*}
0 \to V\cong {\mathfrak V}/t{\mathfrak V}\xrightarrow{i\cong t\times-} 
W\cong {\mathfrak V}/t^{2}{\mathfrak V}
\xrightarrow{p\cong -\otimes_{k[\![t]\!]}k} V\cong {\mathfrak V}/t
{\mathfrak V}\to 0\,.
\end{align*}
Now apply $\Hom_{k[\![t]\!]Q}({\mathfrak V},-)$ to this exact sequence to 
obtain the following long exact sequence of $k[\![t]\!]$--modules, with 
$\delta$ denoting the connecting homomorphism:
\begin{align*}
\xymatrix{
0\ar[r]&
\Hom_{k[\![t]\!]Q}({\mathfrak V},V)\ar[r]&
\Hom_{k[\![t]\!]Q}({\mathfrak V},W)\ar[rr]&&
\Hom_{k[\![t]\!]Q}({\mathfrak V},V)\ar[r]^-{\delta}&{\hphantom  x}\\
\ar[r]&\Ext^{1}_{k[\![t]\!]Q}({\mathfrak V},V)\ar[r]&
\Ext^{1}_{k[\![t]\!]Q}({\mathfrak V},W)\ar[rr]^{\Ext^{1}_{k[\![t]\!]Q}
({\mathfrak V},p)}&&
\Ext^{1}_{k[\![t]\!]Q}({\mathfrak V},V)\ar[r]&
0
}
\end{align*}
The map 
$\pi=\Ext^{1}_{k[\![t]\!]Q}({\mathfrak V},p)$ is the same as the 
projection alluded to above, which we wish to show is an 
isomorphism. Using the various identifications, we may rewrite this long 
exact sequence as
\begin{align*}
\xymatrix{
0\ar[r]&
\End_{Q}(V)\ar[r]&
\End_{({k[\![t]\!]}/{(t^{2})})Q}(W)\ar[r]&
\End_{Q}(V)\ar[r]^-{\delta}&{}\\
\ar[r]&\Ext^{1}_{Q}(V,V)\ar[r]&
\Ext^{1}_{k[\![t]\!]Q}({\mathfrak V},W)\ar[r]^{\pi}&
\Ext^{1}_{Q}(V,V)\ar[r]&
0
}
\end{align*} 
As $\d$ is a Schur root, and $\dim\, T^{1}(V) = \dim\, \Ext^{1}_{Q}(V,V) = 1$, 
we see that $\pi$ is an isomorphism if and only if $\delta\neq 0$ if and 
only if there exists a $Q$--endomorphism of $V$ that cannot be lifted 
to a  $k[\![t]\!]$--linear $Q$--endomorphism of $W$. While these 
considerations apply to any quiver, 
we now show that $\delta\neq 0$, thereby establishing that $\pi$ is indeed an 
isomorphism, for any Schur root of a Dynkin quiver.

By assumption, 
$q(V)=1$
and $\dim\, \Ext^{1}_{Q}(V,V) = 1$, whence $V$ is 
decomposable, say, $V=V'\oplus V''$ for nonzero $Q$--representations $V',V''$. 
It follows from  $\dim\, \End_{Q}(V) =2$ that 
$\End_{Q}(V) \cong \End_{Q}(V')\oplus \End_{Q}(V'')$, and that the 
endomorphism rings of $V',V''$ are one-dimensional, in particular these 
representations are indecomposable. This means that their dimension vectors 
are real Schur roots as well, and so the representations are rigid. From $\Ext^{1}_{Q}(V,V)\cong \Ext^{1}_{Q}(V'\oplus V'',V'\oplus V'')$, 
it then follows that exactly one of the groups $\Ext^{1}_{Q}(V',V'')$ 
or $\Ext^{1}_{Q}(V'',V')$ is nonzero --- and then one-dimensional. 
Assume $\Ext^{1}_{Q}(V',V'')\neq 0$. 
The associated nontrivial extension
\begin{equation}
\label{eq:1}
0\to V'' \xrightarrow{i} W' \xrightarrow{p} V'\to 0
\end{equation}
gives rise to the following nonzero extension class in $\Ext^{1}_{Q}(V,V)$:
\begin{align*}
\begin{matrix}
&&V'&\xrightarrow{\quad =\quad}&V'\\
&&\oplus&&\oplus\\
0&\to&V''&\xrightarrow{\quad i\quad}&W'&\xrightarrow{\quad p\quad}&V'&\to&0\\
&&&&\oplus&&\oplus\\
&&&&V''&\xrightarrow{\quad =\quad}&V''
\end{matrix}
\end{align*}
Note that $W'$ has dimension vector $\d$, as that is the sum of the dimension vectors of $V''$ and $V'$, equal to the dimension vector of $V$. It is now a general fact that $V=V'\oplus V''$ deforms into the middle term $W'$, for any extension. As the sequence does not split, $W'\not\cong V$, and, as $V$ has a onedimensional semi-universal deformation, $W'$ must be the indecomposable representation of dimension vector $\d$. 
Using the observation following (\ref{al:1}), the $k[\![t]\!]$--module structure on the middle term $W=V'\oplus W'\oplus V''$ is as follows:
\begin{align*}
t(V') &= 0\quad,\quad t|_{W'}= p\quad, \quad t|_{V''}= i\,. 
\end{align*}
With $W'$ an indecomposable $Q$--representation and the action of $t$ as described, it follows easily that $W=V''\oplus W'\oplus V'$ is indecomposable as a $Q$--representation over $k[\![t]\!]$. Accordingly, its endomorphism ring $\End_{({k[\![t]\!]}/{(t^{2})})Q}(W)$ contains only the trivial idempotents, thus none of the idempotents in $\End_{Q}(V)$ that corresponds to the projections onto the indecomposable factors of $V$ can be lifted, and the natural ring homomorphism $\End_{({k[\![t]\!]}/{(t^{2})})Q}(W)\to \End_{Q}(V)$ is not surjective. This yields the claim.\eop
 
\begin{corollary} If $Q$ is a Dynkin quiver and $\d$ is a real root 
of $Q$ then the discriminant in $\Rep(Q,\d)$ is a linear
free divisor.
\end{corollary}
\ni\Proof By Gabriel's theorem $Q$ is of finite representation type. 
Therefore at a generic point $V$ on each irreducible component of $D$,
any deformation of $V$ inside $D$ is trivial. Thus $T^1(V)$ is 1-dimensional.
\eop
\ni 
Everything we have said so far only depends on the 
{\em support\/} of the dimension vector $\d$, that is, the full subquiver 
whose nodes are those $x\in Q_{0}$ with $\d(x)\neq 0$. A dimension vector 
is {\em sincere\/} if its support is all of $Q_{0}$.
\section{A Criterion for $D$ to be a Linear Free Divisor}\label{crit}
The group $\Gl(Q,\d)$ acts on the ring $R$ of polynomial functions on 
$\Rep(Q,\d)$ by the contragredient action, as described earlier in 
Section \ref{defs}. 
A polynomial $f\in R$ is a {\it semi-invariant of weight
$\chi$}, where $\chi$ is a character of $\Gl(Q,\d)$, 
if for all $g\in \Gl(Q,\d)$ we have $g\cdot f = \chi(g)f$. As
the characters of $\Gl_n(k)$ are just integral powers of $\det$, 
the characters
of $\Gl(Q,\d)$ are in bijection with elements of $\ZZ^{Q_0}$. The 
{\it weight} $w(f)$ of a semi-invariant $f$ is usually 
identified with the image in $\ZZ^{Q_0}$ of its associated character. 
\begin{theorem}\label{sakim} {\em (Sato-Kimura \cite{saki})} Let the connected
algebraic group $G$ act on the vector space $V$. If there is an open orbit
then the ring $SI(G,V)$ 
spanned by the semi-invariants
is a polynomial ring:
$$SI(G,V)=k[f_1,\ld,f_s]$$
for some collection of algebraically independent and irreducible 
semi-invariants $f_1,\ld,f_s$. Moreover if $f_i\in SI(G,V)_{\chi_i}$
then the $\chi_i$ are linearly independent in the space of characters of
$G$. 
\eop
\end{theorem}
\begin{corollary}
Under the assumptions of the theorem, the set of characters $\chi$ such that
$SI(G,V)_{\chi}\neq 0$ forms a free abelian semigroup, isomorphic to 
$\NN^{s}$. In particular, if $f$ is any semi-invariant, of weight $\chi$, 
then $f=uf_{1}^{a_{1}}\cdots f_{s}^{a_{s}}$, where $u$ is a unit in $k$ 
and the $a_{i}\ge 0$ are the unique integers such that $\chi = 
\sum_{i=1}^{s}a_{i}\chi_{i}$ in the space of characters of $G$.\eop
\end{corollary}
\ni
Suppose that $\d$ is a real Schur root of $Q$, and let $D$ be the discriminant 
in $\Rep(Q,\d)$. As $D$ is preserved under the action of
$\Gl(Q,\d)$, its canonical equation $\Delta$  is a semi-invariant. 
If $V\notin D$, and
$f$ is a non-zero semi-invariant, then $f(V)$ cannot vanish; if
it did, then it would vanish everywhere on the orbit of $V$, which is dense. 
In other words, the zero locus of any semi-invariant must be contained 
in the discriminant. In particular, with the $f_i$ as in \ref{sakim}, 
$f_1\cdots f_s$ is necessarily a reduced equation for $D$, and so $\Delta = 
uf_{1}^{a_{1}}\cdots f_{s}^{a_{s}}$, with $u$ a unit in $k$, and uniquely 
determined integers $a_{i}>0$. 

\ni Moreover, Kac has shown in \cite[p.153]{kac2} that the discriminant 
for a real Schur root $\d$ contains precisely $n-1$ irreducible components, 
where $n$ is the number of nodes in the support of  $\d$, thus, there 
are $s=n-1$ fundamental semi-invariants $f_{i}$ in $SI(\Gl(Q,\d), 
\Rep(Q,\d))$. This gives us a first combinatorial criterion for $D$ to 
be a linear free divisor.
\begin{proposition}
\label{prop:crit}
Suppose that $\d$ is a real Schur root of $Q$, supported on $n$ nodes. 
Assume further that $g_{1},...,g_{n-1}$ are semi-invariants on $\Rep(Q,\d)$ 
with linearly independent weights $w_{i}=w(g_{i})$. If the weight of the 
discriminant $D$ 
satisfies $w(D) = \sum_{i=1}^{n-1}a_{i}w_{i}$, for integers $a_{i}\ge 1$, 
then 
$\Delta = ug_{1}^{a_{1}}\cdots g_{n-1}^{a_{n-1}}$ for some unit $u\in k$.
If we know further that the weights 
$w_{i}$ generate the semigroup of all weights occurring in 
$SI(\Gl(Q,\d), \Rep(Q,\d))$, then the $g_{i}$ constitute the reduced 
equations of the components of $D$, and $D$ is a linear free divisor 
if and only if each $a_{i}=1$.\eop
\end{proposition}
\ni Derksen and Weyman in \cite{derksenweyman} describe in general the 
semigroup of weights occurring in $SI(\Gl(Q,\d), \Rep(Q,\d))$ through 
a single 
equation\footnote{Namely that the ordinary scalar product of the weight 
of a semi-invariant with the dimension vector $\d$ has to vanish, that is, 
$w\cdot\d=0$.}
 and integral inequalities that depend upon the dimension vectors of generic 
subrepresentations, whence the criterion can be applied, at least in principle. 
We may as well turn the criterion around to determine all semi-invarants 
if we already know that $D$ is a linear free divisor, such as for real 
roots whose support is a Dynkin quiver:
\begin{corollary}
Assume the discriminant $D$ in $ \Rep(Q,\d)$, for $\d$ a real Schur root, 
is a free divisor and its canonical equation factors as $\Delta= g_{1}
\cdots g_{n-1}$ for semi-invariant polynomials $g_{i}$ with linearly 
independent weights. If $n$ is the number of nodes in the support of $\d$, 
then the factors $g_{i}$ are algebraically independent and irreducible 
polynomials that generate the ring of semi-invariants $SI(\Gl(Q,\d), 
\Rep(Q,\d))$.\eop
\end{corollary}

\ni Using yet another result of Schofield \cite{schofield}, one may find the 
weights of all semi-invariants --- indeed the semi-invariants themselves, 
as we will discuss in more detail later, see Section \ref{eqns}.
Suppose that $\e$ is a dimension vector such that $\la \e,\d\ra=0$. In the 
exact sequence (\ref{fes}), the matrix $d^W_V$ is now square. 
We define a polynomial function 
$c:\Rep(Q,\e)\times\Rep(Q,\d)\to k$
by $c(W,V)=\det\,d^W_V$. The map 
$$\Rep(Q,\e)\times\Rep(Q,\d)\to \Hom_k\Bigl(\prod_{x\in Q_0}
\Hom(k^{\e(x)},k^{\d(x)}),
\prod_{\pp\in Q_1}\Hom_k(k^{\e(t\pp)},k^{\d(h\pp)}\Bigr)$$
sending $(W,V)$ to $d^W_V$ is $\Gl(\e)\times Gl(\d)$-equivariant, 
and it follows that for fixed
$W$, the map $c^W:=c(W,\ )$
represents a semi-invariant polynomial on $\Rep(Q,\d)$.
\begin{theorem}{\em(\cite{schofield} 4.3)}
\label{semiinvs}
Let $Q$ be a quiver without oriented cycles, and
let $\d$ be a sincere real Schur root for $Q$. The polynomials $c^W$ with
$\la W,V\ra=0$ span the ring of semi-invariants $SI(\Gl(\d),\Rep(Q,\d)$.\eop
\end{theorem}

\ni The weights of these semi-invariants, as well as that of the discriminant, 
are then easily established, using essentially the same argument as in 
\cite[1.4]{schofield}. To formulate it succinctly, we introduce the 
{\em in--degree\/} $\indeg_{\d}$ and the {\em out--degree\/} $\outdeg_{\d}$ 
of $\d$ as the dimension vectors
\begin{align}
\label{al:inout}
\indeg_{\d}(x) &=\sum_{ \pp\in Q_{1} : h\pp=x} \d(t\pp) \quad,\quad
\outdeg_{\d}(x) = \sum_{\pp\in Q_{1}: t\pp=x} \d(h\pp)\quad,\quad\text
{for $x\in Q_{0}$.}
\end{align}
In terms of the Euler matrix $E$ of $Q$ (see Section \ref{background}) one 
has 
\begin{align*}
\indeg_{\d}  = \d - \d E\quad,\quad \outdeg_{\d} = \d - \d E^{T}\,.
\end{align*}

\begin{lemma}
\label{lem:weights}
Let $\d, \mathbf e$ be dimension vectors for the quiver $Q$ with 
$\langle\e,\d\rangle =0$. The weight of the $\Gl(\e)\times\Gl(\d)$ 
semi-invariant polynomial $c(W,V)$ in the character group $\bbbz^{Q_{0}}
\times \bbbz^{Q_{0}}$ is 
\begin{align*}
w(c(W,V))&=(\d - \outdeg_{\d}, -\e +\indeg_{\e}) = (\d E^{T},-\e E)\,, 
\intertext{while that of the $\Gl(\d)$ semi-invariant $c^{W}$ in 
$\bbbz^{Q_{0}}$ is}
w(c^{W}) &= -\e +\indeg_{\e} = -\e E
\intertext{and the weight of the discriminant in $\Rep(Q,\d)$ equals}
w(\Delta) &= \indeg_{\d} -\outdeg_{\d} = \d(E^{T}-E)\,.
\end{align*}
%
\end{lemma}
\ni\proof
Let $V,W$ be two representations with dimension vectors $\d,\mathbf e$ such 
that $\langle\e,\d\rangle =0$. The map $d^{W}_{V}$ can be viewed as a 
linear map
\begin{align*}
d^{W}_{V}:\bigoplus_{x\in Q_{0}}V_{x}\otimes W_{x}^{*} &\too 
\bigoplus_{\pp\in Q_{1}}V_{h\pp}\otimes W_{t\pp}^{*}\,,
\end{align*}
where $(-)^{*}$ denotes the $k$--dual. Denoting by $\Lambda(-)$ the highest 
exterior power of a vector space, and observing that 
\begin{align*}
\Lambda(U^{*})&\cong \Lambda(U)^{*}\,,\,\Lambda(U\oplus U')\cong 
\Lambda(U)\otimes \Lambda(U')\,,\,
\Lambda(U\otimes U') \cong \Lambda(U)^{\dim U'}\otimes  \Lambda(U')^{\dim 
U}\,,
\end{align*}
for vector spaces $U,U'$,
the determinant of $d^{W}_{V}$ can be represented as
\begin{align*}
\det d^{W}_{V}\cong \Lambda (d^{W}_{V})&: \bigotimes_{x\in Q_{0}}
\Lambda(V_{x})^{\e(x)}\otimes\Lambda(W_{x}^{*})^{\d(x)}\too
 \bigotimes_{\pp\in Q_{1}}\Lambda(V_{h\pp})^{\e(t\pp)}\otimes
\Lambda(W_{t\pp}^{*})^{\d(h\pp)}\,.
\end{align*}
One reads off that as a semi-invariant for $\Gl(\e)\times\Gl(\d)$ the 
determinant of $d^{W}_{V}$ transforms according to
\begin{align*}
\left(\prod_{\pp\in Q_{1}}\det\bigl(\Gl(\d(h\pp))\bigr)^{\e(t\pp)}
\det\bigl(\Gl(\e(t\pp))\bigr)^{-\d(h\pp)}\right)
\left(\prod_{x\in Q_{0}}\det\bigl(\Gl(\d(x))\bigr)^
{-\e(x)}\det\bigl(\Gl(\e(x))\bigr)^{\d(x)}\right)
\end{align*}
thus, its weight, in the character group $\bbbz^{Q_{0}}\times\bbbz^{Q_{0}}$ 
of $\Gl(\e)\times\Gl(\d)$, is given on a pair of nodes $(y,x)$ by
\begin{align*}
w(\det d^{W}_{V})(y,x) &= \d(y) -\sum_{t\pp=y} \d(h\pp) -\e(x) + 
\sum_{h\pp=x} {\e}(t\pp)  
\intertext{thus,}
w(\det d^{W}_{V})&= (\d -\outdeg_{\d}, -\e + \indeg_{\e}) = (\d E^{T}, -\e E)
\quad\in\quad \bbbz^{Q_{0}\times Q_{0}}.
\end{align*}
For $V=W$, the diagonal summand $k\subseteq \oplus_{x\in Q_{0}}\Hom(V_{x},
V_{x})$ does not contribute to the weight of the determinant, and restricting 
$w(\det d^{V}_{V})$ to the diagonal $y=x$ yields the claimed formula 
for the discriminant.
\eop

Now we are ready to study some examples.
\section{Examples} 
To illustrate the results and to exhibit explicit linear free divisors 
arising from Dynkin quivers, we concentrate mainly on the most complicated 
ones, those corresponding to the {\em highest root\/} of a Dynkin diagram 
viewed as the dimension vector of some Dynkin quiver. Recall that the 
connected Dynkin diagrams are in natural bijection with the {\em binary 
polyhedral groups\/}, the congruence classes of finite subgroups of 
$\Sl(2,\CC)$. One has the following simple relation between the dimension 
of the representation variety associated to the highest root and the 
order of the corresponding finite group.

\begin{lemma}
Let $Q$ be a connected Dynkin quiver, $\d$ the highest root of the underlying 
Dynkin diagram, and $\Gamma$ the associated binary polyhedral group. 
The dimension of $\Rep(Q,\d)$, equal to the degree of the discriminant $D$, 
is then
$
\dim\,\Rep(Q,\d) = |\Gamma| - 2\,.
$
\end{lemma}

\ni\proof By the McKay correspondence, the components $\d(x)$ of the highest 
root are in bijection with the dimensions of the isomorphism classes of 
irreducible and nontrivial representations of $\Gamma$. Accordingly,
\begin{align*}
|\Gamma| = 1+\sum_{x\in Q_{0}}\d(x)^{2} = 2 + \dim\, \pgl(\d) = 2 + 
\dim\,\Rep(Q,\d)\,.
\end{align*}
\eop

\begin{example}
{\em Let $Q$ be a Dynkin quiver of type $A_n$ with any orientation, 
and let $\d$ be its highest root, the dimension vector assigning 1 at 
each vertex. Then 
$\Rep(Q,\d)$ can be identified with $k^{|Q_1|}=k^{n-1}$ by associating 
to each morphism
its $1\times 1$ matrix. Each of the coordinates is a semi-invariant, 
and
$D$ is the normal crossing divisor in $n-1$ variables. Notice that $D$ 
is independent of
the orientation of the arrows.}
\end{example}
\begin{example}
{\em
Consider the two Dynkin quivers $Q^{(1)}$ and $Q^{(2)}$ 
of type $E_6$ with the highest root as 
dimension vector as
shown. Each space $\Rep(Q^{(i)},\d)$ has 
dimension $22 = 24-2$, as the corresponding 
binary tetrahedral group has order 24. 

$$
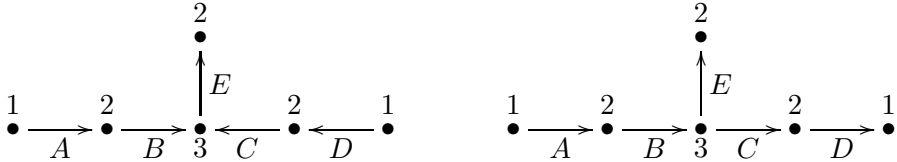
\xymatrix{&&\bullet\save*++!D{2}\restore&&\\
\bullet\save*++!D{1}\restore
\ar[r]_{\displaystyle A}&\bullet\save*++!D{2}\restore\ar[r]_{\displaystyle B}&
\bullet\save*+!U{3}\restore\ar[u]_{\displaystyle E}
&
\bullet\save*++!D{2}\restore
\ar[l]^{\displaystyle C}&
\bullet\save*++!D{1}\restore\ar[l]^{\displaystyle D}
}
\hskip 30pt
\xymatrix{&&\bullet\save*++!D{2}\restore&&\\
\bullet\save*++!D{1}\restore
\ar[r]_{\displaystyle A}&\bullet\save*++!D{2}\restore\ar[r]_{\displaystyle B}&
\bullet\save*+!U{3}\restore\ar[u]_{\displaystyle E}\ar[r]_{\displaystyle C}
&
\bullet\save*++!D{2}\restore
\ar[r]_{\displaystyle D}&
\bullet\save*++!D{1}\restore
}
$$
\vskip 10pt
One sees easily that codimension 1 degeneracies are
given, for $Q^{(1)}$, by the vanishing of any of\footnote{We indicate by 
$X|Y$ the concatenation of two matrices $X,Y$ with the same number of rows.}
$$\det[EB],\ \det[EC],\ \det[B|CD],\ \det[BA|C],\ \det[EBA|ECD]\,.$$
The third of these measures the independence of the images of $B$ and $CD$
in the 3-dimensional space attached to the central node; the fourth and fifth
are to be understood similarly. The degrees of the corresponding equations, 
equal to $4, 4, 4, 4$, and $6$, 
add to 22, and their weights are easily seen to be linearly independent. 
Thus these form a complete list of the factors, and
the linear free divisor $D$ is the union of these
five, necessarily irreducible components.

For $Q^{(2)}$, four codimension 1 degeneracies are defined by
the vanishing of 
$$\det[EB],\det[CB], 
\det\left[\begin{array}{c}E\\DC\end{array}\right], DCBA\,.$$ 
One further degeneracy
is easier to describe verbally than by an equation 
(however, see Section \ref{eqns} and in particular Example
\ref{exe7} below): it is the failure
of general position,
in the 3-dimensional space at the central node, of the three lines 
$\im(BA),\ker(E),\ker(C)$. 

In both cases, each equation of degree 4 has 12 monomials and the equation
of degree 6 has 48. Moreover, the complements of the discriminants
$D^{(1)}\subset\Rep(Q^{(1)},\d)$ and $D^{(2)}\subset \Rep(Q^{(2)},\d)$
are isomorphic to one another, being orbits, with trivial isotropy, 
of the groups $\PP\mbox{Gl}(Q^{(i)},\d)$, which are 
themselves isomorphic to one another. 
However, the two discriminants are not isomorphic. 
Essentially, this is because the equations involve 
different numbers of variables.
In the
first case, the five equations involve, respectively, 12, 12, 14, 14, and 22
variables, while in the second the five equations involve 12, 12, 14, 16 and
20 variables. A {\it Macaulay} calculation confirms that the spaces
of vector fields with constant coefficients
tangent to the germs at $0\in\Rep(Q^{(i)},\d)$ of the five components have 
dimensions 10, 10, 8, 8, and 0 in the first case, and 10, 10, 8, 6 and 2
in the second. Any isomorphism $D^{(1)}\cong D^{(2)}$ must map $0$ to $0$,
since because of the presence of the Euler field, 
in each case $0$ is the only point where all of the vector fields 
in $\Der(-\log D^{(i)})$ vanish. It follows that these dimensions are 
geometrical 
invariants: the dimension corresponding to the irreducible component 
$D^{(i)}_j$ of $D^{(i)}$ is
the maximum dimension of a non-singular factor in a product decomposition
$(D^{(i)}_j,0)\cong (E^{(i)}_j,0)\times (F^{(i)}_j,0)$. 
} 
\end{example}
\begin{proposition}\label{nn+1} 
Let $Q$ be the quiver whose nodes consist of $n+1$ sources surrounding 
one sink, with an arrow going from each source to the sink. The 
discriminant with respect to the dimension vector that assigns $1$ to 
each of the sources and $n$ to the sink is a linear free divisor. It 
is of the form $\Delta_{1}\cdots \Delta_{n+1}$, where the $\Delta_{i}$ 
are the maximal minors of a generic $n\times (n+1)$--matrix.
\end{proposition}
\ni\Proof We can identify $\Rep(Q,\d)$ with the space of $n\times(n+1)$--matrices, with the matrix of each of the arrows determining a column. 
The degree of the discriminant $D$ equals $n(n+1)$. The generic representation 
describes $n+1$ distinct lines in a vector space  of dimension $n$, with 
no $n$ of them lying in a hyperplane. Such a representation is indecomposable 
and lies in an open orbit, with the group $\Gl(n)$ acting transitively 
on the set of such line configurations in general position. Accordingly, 
the dimension vector is a real Schur root.
There are $n+1$ codimension 1 degeneracies, each one determined by the
vanishing of an $n\times n$ minor of the $n\times (n+1)$ matrix. The product
of these minors has degree $n(n+1)$, equal to the degree of $D$, and the 
weights, assigning $-1$ to each source contributing to the minor, $0$ 
to the remaining source, and $1$ to the sink, are clearly linearly 
independent. Thus each
is present in $\det\,\tilde d^M_M$ with multiplicity 1.\eop 
Note that from Theorem \ref{sakim} we 
recover the classical result that these maximal minors are algebraically 
independent.

\begin{example}{\em
Consider the four
quivers shown below, 
in which the underlying 
undirected graph is the extended Dynkin diagram of type 
$\widetilde D_{4}$. 
Assign to each the dimension
vector with 1 at each outer node and 3 at the central node. 
According to Kac's result quoted as Proposition \ref{kacs} above, 
the dimension vector shown is a real root. In (i)--(iii), it is also a 
Schur root, but in case (iv), it is not. In case (i),
the discriminant is a linear free divisor, according to 
Proposition \ref{nn+1} above, but in cases (ii) and (iii) this fails. 
In case (iv), the discriminant is the whole space, and there is 
no rigid representation.
$$
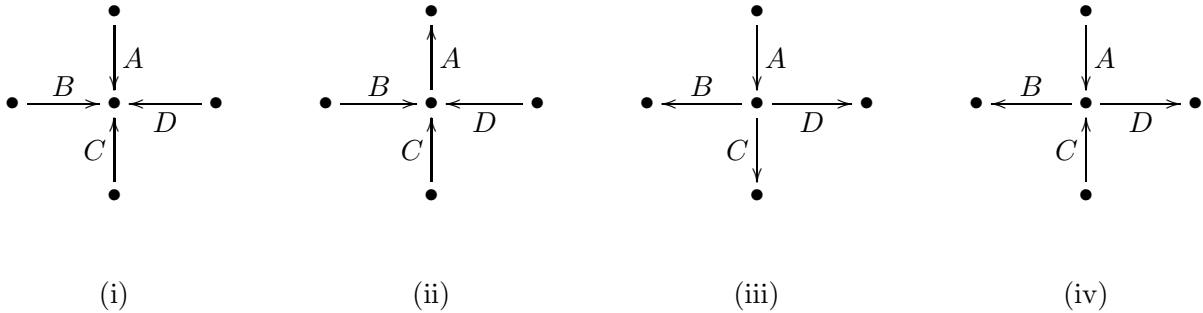
\xymatrix{&\bullet\ar[d]^{\displaystyle A}&\\
\bullet\ar[r]^{\displaystyle B}&\bullet&\bullet\ar[l]^{\ds D}\\
&\bullet\ar[u]^{\ds C}\\
&\txt{(i)}
}
\hskip 30pt
\xymatrix{&\bullet&\\
\bullet\ar[r]^{\displaystyle B}&\bullet\ar[u]_{\ds A}&\bullet\ar[l]^{\ds D}\\
&\bullet\ar[u]^{\ds C}\\
&\txt{(ii)}
}
\hskip 30pt
\xymatrix{&\bullet\ar[d]^{\displaystyle A}&\\
\bullet&\bullet\ar[l]_{\ds B}\ar[r]_{\ds D}\ar[d]_{\ds C}&\bullet\\
&\bullet\\
&\txt{(iii)}
}
\hskip 30pt
\xymatrix{&\bullet\ar[d]^{\displaystyle A}&\\
\bullet&\bullet\ar[l]_{\ds B}\ar[r]_{\ds D}&\bullet\\
&\bullet\ar[u]^{\ds C}\\
&\txt{(iv)}
}
$$
In case (ii), there is a modulus attached to the codimension 1 
degeneracy in which the images of $B,C$ and $D$ lie in a plane $P$; these 
three
lines, together with the fourth line $P\cap\ker\, A$, determine
a cross-ratio. Any representation $V$ of this type therefore has $T^1_V$ of
dimension (at least) 2, and so the multiplicity of the corresponding
component in $D$ is also at least 2. In fact it is exactly 2: the remaining
three components of $D$ are $\det\,AB,\det\,AC,\det\,AD$, each of degree 2.
Together with twice the degree of $\det[B|C|D]$ these add up to 12, the 
degree of the (non-reduced) equation $\det\,\tilde d^M_M$
of $D$. As the four components described have linearly independent weights, 
the multiplicity of the non-reduced component is exactly 2.

Case (iii), obtained by reversing all of the arrows, is dual to (ii): here 
the non-reduced component of $D$ is where the kernels of the three
outgoing arrows $B,C,D$ meet along a line $L$. Together with the plane
$L+\im\,A$, these make four planes in the pencil of planes containing $L$,
and thus once again determine a cross ratio. 

In the fourth quiver, the given dimension vector is not a Schur root. For
there is no open orbit. In a general representation $V$,
$\im\ A$ and $\im\ C$ span a plane $P$. The intersections with $P$ of 
$\ker\,B$ and $\ker\,D$ determine two further lines in $P$, and thus a 
cross-ratio. Since thus $\dim\,\Ext^1_Q(V,V)\geq 1$, it follows that 
$\dim\,\Hom_Q(V,V)\geq 2$, and $V$ must be decomposable. Indeed, it is 
easily verified that the intersection $\ker\,B\cap \ker\,D$, concentrated 
on the central node, splits off.
By Kac's theorem, there is exactly one orbit of indecomposable 
representations. We invite the reader to find it.  
} 
\end{example}

\begin{proposition} Suppose that $\d$ is a real Schur root of the 
connected quiver $Q$, and let $Q^{\opp}$ be obtained from $Q$ by reversing
all of the arrows. If the discriminant in 
$\Rep(Q,\d)$ is a linear free divisor then the same holds in 
$\Rep(Q^{\opp},\d)$.
\end{proposition}
\ni\Proof This is essentially projective duality. Transposition 
determines an isomorphism of representation spaces 
$\Rep(Q,\d) \to \Rep(Q^{\opp},\d)$ which maps
orbits to orbits. \eop

\begin{example}{\em Suppose $Q$ is a quiver and $x\in Q_0$ is a node.
Construct a new
quiver $Q_x$ by replacing the node $x$ by a pair of nodes $x',x''$ 
connected by an arrow $F$ from $x'$ to $x''$,
and attaching the arrows previously attached to $x$ either to $x'$ or to $x''$.
Two possible outcomes of this process are shown in the figure below.
If $\d$ is any dimension vector for $Q$, we define a dimension vector $\d_x$
for $Q_x$ by setting $\d_x(y)=\d(y)$ if $y\neq x',x''$, 
$\d_x(x')=\d_x(x'')=\d(x)$. 
Then $\la \d_x,\d_x\ra=\la\d,\d\ra$. If the generic representation 
in $\Rep(Q,\d)$ is indecomposable, then the same is true in $\Rep(Q_x,\d_x)$,
since generically $V(F)$ is an isomorphism. So it is reasonable to 
hope that if $\la\d,\d\ra=1$ and $D\subset\Rep(Q,\d)$ is a linear free divisor, 
then  
the discriminant in $\Rep(Q_x,\d_x)$ is also a linear free divisor. 
The following
examples show that this is sometimes but not always the case.\\

\ni The quivers $Q_2$ and $Q_3$ shown below are obtained from $Q_1$
by the operation just described. Assign to $Q_1$ the dimension vector
$\d$ with 1's at all the sources and 4 at the central sink, and define $\d_x$
accordingly. By \ref{nn+1}, the discriminant $D_1\subset\Rep(Q_1,\d)$ 
is a linear free divisor with
components given by the vanishing of
$$
\det[A|B|C|D], \det[A|B|C|E],\det[A|B|D|E], \det[A|C|D|E], \det[B|C|D|E]\,.
$$
In $\Rep(Q_2,\d_x)$, these become
$$
\det[FA|FB|FC|D], \det[FA|FB|FC|E], \det[FA|FB|D|E], \det[FA|FC|D|E],
\det[FB|FC|D|E], \det\,F\,.
$$
In $\Rep(Q_3,\d_{x})$, they become
$$\det[A|B|C|FD], \det[A|B|C|FE], \det[A|B|FD|FE], \det[A|C|FD|FE], 
\det[B|C|FD|FE], \det\,F\,.$$
The degrees of the (reduced) discriminants $D_2\subset
\Rep(Q_2,\d_x)$ and $D_3\subset\Rep(Q_3,\d_x)$ are thus 
$36$ and $32$ respectively. So $D_2$ is a linear free divisor,
whereas $D_3$ is not. The exponent of $\det\,F$ in the canonical equation 
$\Delta_3$ is 2.
$$
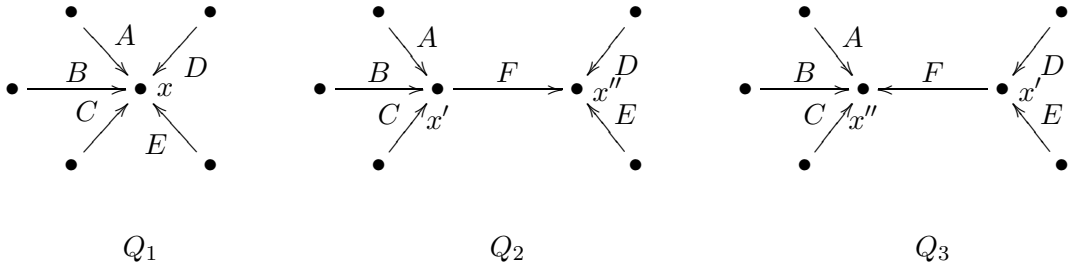
\xymatrix@R=0.25in@C=0.15in{&\bullet\ar[dr]^{\ds A}&&\bullet\ar[dl]^{\ds D}\\
\bullet\ar[rr]^{\ds B}&&\bullet
\save*++!L
{\ds x}
\restore\\
&\bullet\ar[ur]^{\ds C}&&\bullet\ar[ul]^{\ds E}
\\
&&\txt{$Q_1$}}
\hskip 30pt
\xymatrix@R=0.25in@C=0.15in{&\bullet\ar[dr]^{\ds A}&&&&\bullet\ar[dl]^{\ds D}\\
\bullet\ar[rr]^{\ds B}&&\bullet
\save*++!U
{\ds x'}\restore
\ar[rr]^{\ds F}
&&\bullet
\save*++!L{\ds x''}\restore
\\
&\bullet\ar[ur]^{\ds C}&&&&\bullet\ar[ul]_{\ds E}
\\
&&&\txt{$Q_2$}
}
\hskip 30pt
\xymatrix@R=0.25in@C=0.15in{&\bullet\ar[dr]^{\ds A}&&&&\bullet\ar[dl]^{\ds D}\\
\bullet\ar[rr]^{\ds B}&&\bullet
\save*++!U
{\ds x''}\restore
&&\bullet
\save*++!L{\ds x'}\restore\ar[ll]_{\ds F}
\\
&\bullet\ar[ur]^{\ds C}&&&&\bullet\ar[ul]_{\ds E}
\\
&&&\txt{$Q_3$}
}
$$
One can easily show, by the same technique of counting degrees, that if one 
performs this operation
on the
central node in the quiver of Proposition \ref{nn+1}, then one obtains a linear
free divisor if and only if just two of the arrows coming from the outer nodes
are attached to $x''$, and the rest are attached to $x'$.

By applying the same construction to Dynkin quivers and their roots, one
can obtain further examples of linear free divisors. In particular, one easily deals with the case $D_{n}$ in this way:}
\end{example}

\begin{proposition}
Let $Q$ be the Dynkin quiver of type $D_{n}$ with the following orientation 
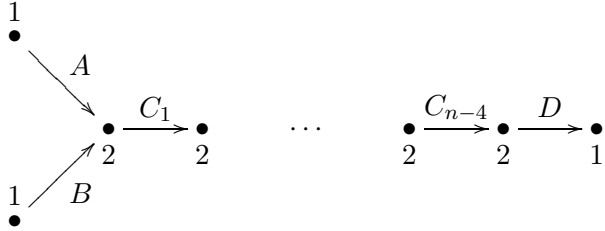
\begin{align*}
\xymatrix{
\bullet\save*++!D{ 1} \restore
\ar[dr]^{\displaystyle A}\\
&
\bullet \save*++!U{ 2} \restore
\ar[r]^{\displaystyle C_{1}} &
\bullet \save*++!U{ 2} \restore&\cdots&
\bullet \save*++!U{2} \restore
\ar[r]^{\displaystyle C_{n-4}} &
\bullet \save*++!U{ 2} \restore
\ar[r]^{\displaystyle D_{\vphantom 2}} &
\bullet \save*++!U{1} \restore\\
\bullet\save*++!D{ 1} \restore
\ar[ur]_{\displaystyle B}\\
}
\end{align*}
The indicated dimension vector $\d$ is the highest root of $\overline D_{n}$. 
The discriminant in $\Rep(Q,\d)$ is a linear free divisor of degree $4n-10$ 
with $n-1$ factors 
\begin{align*}
\det[A|B], \det C_{1},\ldots, \det C_{n-4} , DC_{n-4}\cdots C_{1}A, 
DC_{n-4}\cdots C_{1}B\,,
\end{align*}
where the first $n-3$ factors are of degree $2$, the last two of degree $n-2$.

Changing the orientation of arrows in $Q$ results in an isomorphic linear 
free divisor.
\end{proposition}

\ni \proof The criterion \ref{prop:crit} shows immediately that the factors 
are correct, 
as they represent semi-invariants with linearly independent weights. 
For the last assertion, note that changing the direction of the arrow 
underlying the matrix $C_{i}$, say, results in the same linear free divisor 
as the one already established, provided one replaces $C_{i}$ by its 
adjoint matrix. Similarly, changing, say, the direction of the arrow 
underlying $A$, amounts to replacing $A = (a_{1},a_{2})$ by $A' = 
(a_{2},-a_{1})$ in the above factors, and the situation for $B,D$ is 
analogous.
\section{Equations for $D$}\label{eqns}
To find equations for $D$ in general, one can use the following recipe 
due to Schofield \cite{schofield} that is based on his result \ref{semiinvs} 
above.
We quote it in the slightly simplified form that is all that we require 
here. Assume that $Q$ is a finite connected quiver without oriented cycles 
and fix the sincere real Schur root $\d$ and a generic representation 
$V\in \Rep(Q,\d)$.

To apply \ref{semiinvs}, one looks for roots $\e$ of $Q$ such that 
$\la \e,\d\ra=0$, and computes, for generic $W$ in $\Rep(Q,\e)$, the 
polynomial $c^W$. If $\Hom_{Q}(W,V)\neq 0$, then the square matrix 
underlying $c^{W}_{V}$ has a nontrivial kernel and $c^{W}$ vanishes on 
the open orbit, thus, identically. In view of this, one needs only to 
consider representations $W$ that lie in the {\em left\footnote{One may 
as well work throughout with the right orthogonal category $V^{\perp}$, 
the treatment is symmetric.} orthogonal\/} category ${}^{\perp}V$, the 
full subcategory of all those finite dimensional representations $W$ of 
$Q$ that satisfy
\begin{align*}
\Hom_{Q}(W,V)=\Ext^{1}_{Q}(W,V) = 0\,.
\end{align*}
Schofield shows that this left orthogonal category is equivalent to the 
category of finite dimensional representations of some new quiver $Q'$ 
that has $n-1$ nodes and contains no oriented cycles. In \cite{derksenweyman} 
(Lemma 1) it is pointed out that a short 
exact sequence
$$0\to W'\to W\to W''\to 0$$
of representations of $Q$ leads either to the factorisation 
$$c^W=c^{W''}c^{W'}$$
if $\la W',V\ra=\la W'',V\ra = 0$, or to the conclusion that 
$c^W=0$ if $\la W',V\ra <0$. So, if the generic representation in 
$\Rep(Q,\e)$ is not simple in ${}^{\perp}V$, the semi-invariant we 
obtain will either
be zero or a non-trivial product of others. Accordingly, one needs to 
consider only the
$n-1$ simple objects $W$ in  ${}^{\perp}V$, and those must provide the 
factors of the
discriminant via the associated determinants $c^{W}$. 
Indeed, the dimension vectors $\e_{i}$ of the simple objects $W_{i}$, for 
$i=1,...,n-1$ 
form the unique basis of the free abelian semigroup of dimension vectors 
$\NN^{Q'_{0}}$ for ${}^{\perp}V$, and their associated characters 
$\langle \e_{i},?\rangle = w(c^{W_{i}}) = -\e_{i} +\indeg_{\e_{i}} = 
-\e_{i} E$, 
see \ref{lem:weights}, form the unique basis of the free abelian semigroup 
of 
weights for the semi-invariants of $\Rep(Q,\d)$. Conversely, knowing the 
weights 
$w_{i}$ of the generating semi-invariants, one may calculate the dimension 
vectors $\e_{i}$ through $\e_{i} = -w_{i}(E^{-1})$, with $E^{-1}$ as exhibited 
in \ref{lem:einv}.

The map $\NN^{Q'_{0}}\to \NN^{Q_{0}}$ that maps the i$^{th}$ basis 
vector to $\e_{i}$ is an isometry with respect to the Euler forms on $Q'$ 
and $Q$, 
and as the simple representations for $Q'$ have real Schur roots as their 
dimension vectors, the same must 
hold true for the dimension vectors $\e_{i}$. Thus, in case of a Dynkin 
quiver $Q$, 
we simply need to go through the list of positive roots that are 
perpendicular to $\d$ 
and find among them the uniquely determined basis for the semigroup 
$\NN^{Q'_{0}}$. 

More generally, if $\d$ is the dimension vector of a 
{\em preprojective 
or pre-injective\/} representation, as is the case for any Schur root 
of a Dynkin 
quiver, (see, e.g. \cite[VIII.1]{ARS} for the definitions and result), 
then one 
can read off the roots $\e_{i}$ from the Auslander--Reiten quiver of 
$Q$, 
as explained in \cite[Proof of Proposition 2.1]{HU}. In that case, the 
quiver 
$Q'$ is obtained from $Q$ by deletion of a node along with its incident 
arrows 
and possibly some changes in the orientation of the remaining arrows. 
It is noteworthy that conversely for any quiver, any dimension vector 
of a 
preprojective or pre-injective representation is a real Schur root, 
thus 
providing a huge reservoir for potentially linear free divisors. 
Given that in this situation one can easily determine the simple objects 
of the orthogonal category from the Auslander--Reiten quiver, it seems 
reasonable to expect that one should be able to decide in general which 
of these roots give rise to linear free divisors.

We now turn to the two most complex Dynkin quivers, those of type $E_{7}$ 
and $E_{8}$,
and demonstrate how the algorithm described here works in practice.
\begin{example}\label{exe7}
{\em Consider the Dynkin quiver of type $E_7$ with Schur root
as shown - the highest root of $E_7$.
$$
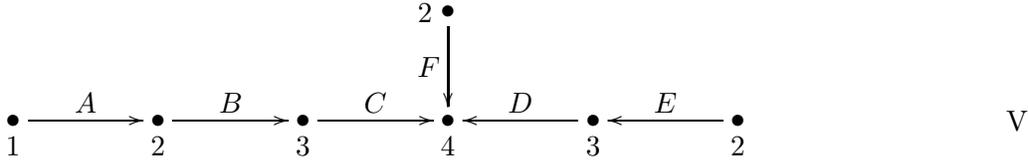
\xymatrix@R=0.4in@C=0.6in{&&&\bullet\save*++!R{2}\restore\ar[d]_{\ds F}\\
\bullet
\save*++!U{1}\restore
\ar[r]^{\ds A}
&
\bullet
\save*++!U{2}\restore
\ar[r]^{\ds B}
&
\bullet
\save*++!U{3}\restore
\ar[r]^{\ds C}
&
\bullet
\save*++!U{4}\restore
&
\bullet
\save*++!U{3}\restore
\ar[l]_{\ds D}
&
\bullet
\save*++!U{2}\restore
\ar[l]_{\ds E}&&\txt{V}
}
$$
\ni The representation space has dimension $46=48-2$ as the associated 
binary polyhedral group, the binary octahedral group, is a double cover 
of the symmetric group on four letters.

By \cite[p.153]{kac2} the discriminant $D$ has 6 irreducible components. Of these,
five may be found by inspection: they are the four described by the equations
$$\det[CBA|D],\ \det[CB|DE],\ \det[F|DE],\ \det[CB|F],$$
and the component corresponding to the degeneracy
$\im\,C\cap\im\,D\cap\im\,F\neq 0$, for which an equation is less obvious. 
One further component remains to be found.
We obtain all of them using Schofield's recipe. Consider first 
$$
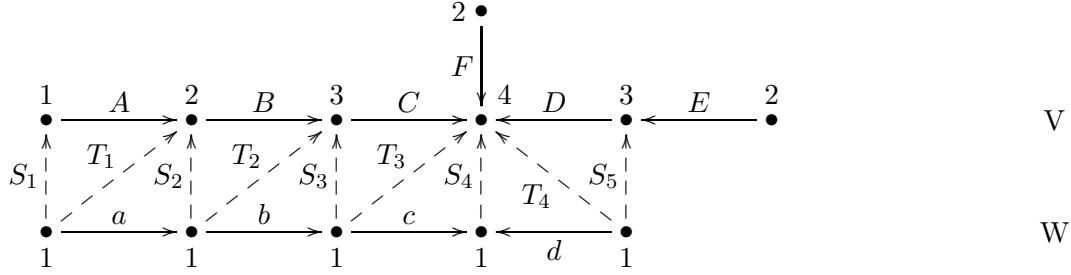
\xymatrix@R=0.4in@C=0.6in{&&&\bullet\save*++!R{2}\restore\ar[d]_{\ds F}\\
\bullet
\save*++!D{1}\restore
\ar[r]^{\ds A}
&
\bullet
\save*++!D{2}\restore
\ar[r]^{\ds B}
&
\bullet
\save*++!D{3}\restore
\ar[r]^{\ds C}
&
\bullet
\save*++!LD{4}\restore
&
\bullet
\save*++!D{3}\restore
\ar[l]_{\ds D}
&
\bullet
\save*++!D{2}\restore
\ar[l]_{\ds E}&&\txt{V}
\\
\bullet
\save*++!U{1}\restore
\ar@{-->}[u]^{\ds S_1}
\ar[r]^{\ds a}
\ar@{-->}[ur]^{\ds T_1}
&
\bullet
\save*++!U{1}\restore
\ar@{-->}[u]^{\ds S_2}
\ar[r]^{\ds b}
\ar@{-->}[ur]^{\ds T_2}
&
\bullet
\save*++!U{1}\restore
\ar@{-->}[u]^{\ds S_3}
\ar[r]^{\ds c}
\ar@{-->}[ur]^{\ds T_3}
&
\bullet
\save*++!U{1}\restore
\ar@{-->}[u]^{\ds S_4}
&
\bullet
\save*++!U{1}\restore
\ar@{-->}[u]^{\ds S_5}
\ar[l]^{\ds d}
\ar@{-->}[ul]^{\ds T_4}
&&&\txt{W}
}
$$
\ni where solid arrows indicate maps within a quiver, and dotted arrows 
indicate 
maps from the quiver $W$ to the quiver $V$ (a convention we adhere to from
now on). Note that the dimension vector $\e$ of $W$ is a root with 
support a Dynkin diagram of type $A_{5}$, the ``type'' of $\e$,  that 
satisfies $\langle \e,\d\rangle =0$.
We have $$d^W_V(S_1,\ld,S_5)=(AS_1-S_2a,BS_2-TS_3b,CS_3-S_4c,DS_5-S_4d).$$
Thus $d^W_V$ has matrix
$$\begin{array}{|c|c|c|c|c|}
\hline
A&-aI_2&0&0&0\\
\hline
0&B&-bI_3&0&0\\
\hline
0&0&C&-cI_4&0\\
\hline
0&0&0&-dI_4&D\\
\hline
\end{array}
$$
where the five columns refer to the five maps $S_1,\ld,S_5$ and the four
rows to the four maps $T_1,\ld,T_4$. Here for each $p,q$ 
we have ordered the natural basis vectors
$E_{ij},{1\leq i\leq q,1\leq j\leq p}$,  of $\Hom(k^p,k^q)$ lexicographically.
Assuming $abc\neq 0$,
row operations transform this 
successively to 
$$\begin{array}{|c|c|c|c|c|}
\hline
A&-aI_2&0&0&0\\
\hline
\frac{1}{a}BA&0&-bI_3&0&0\\
\hline
0&0&C&-cI_4&0\\
\hline
0&0&0&-dI_4&D\\
\hline
\end{array},
\quad
\begin{array}{|c|c|c|c|c|}
\hline
A&-aI_2&0&0&0\\
\hline
\frac{1}{a}BA&0&-bI_3&0&0\\
\hline
\frac{1}{ab}CBA&0&0&-cI_4&0\\
\hline
0&0&0&-dI_4&D\\
\hline
\end{array},
\quad
\begin{array}{|c|c|c|c|c|}
\hline
A&-aI_2&0&0&0\\
\hline
\frac{1}{a}BA&0&-bI_3&0&0\\
\hline
\frac{1}{ab}CBA&0&0&-cI_4&0\\
\hline
\frac{-d}{abc}CBA&0&0&0&D\\
\hline
\end{array}
$$
so that $C(V,W)=\pm d\det[CBA|D]$, and fixing $d\neq 0$ we obtain 
the first of the degeneracies listed
above. Note that the indicated root $\e$ underlying $W$ predicts, by 
\ref{lem:weights}, the following weight of the semi-invariant $c^{W}$:
\begin{align*}
w(c^{W}) &= -\e +\indeg_{\e}:\quad
\begin{matrix}
&&&0\\
-1&0&0&1&-1&0
\end{matrix}
\end{align*}
 which is indeed the weight of $\det[CBA|D]$.
The reader will have no difficulty checking that the next three 
semi-invariants
listed above can be obtained, by the same procedure, from the first three
roots in the diagram
$$
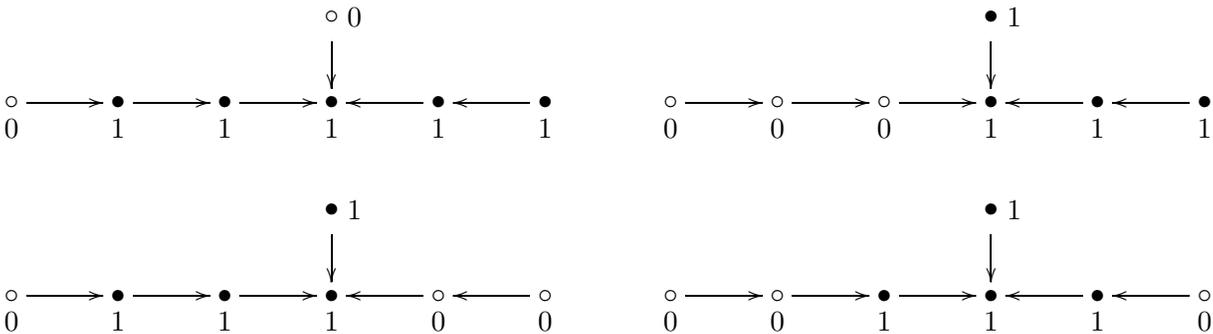
\xymatrix@R=0.3in@C=0.4in{&&&{\circ}\save*++!L{0}\ar[d]\restore\\
{\circ}
\save*++!U{0}\restore
\ar[r]
&
\bullet
\save*++!U{1}\restore
\ar[r]
&
\bullet
\save*++!U{1}\restore
\ar[r]
&
\bullet
\save*++!U{1}\restore
&
\bullet
\save*++!U{1}\restore
\ar[l]
&
\bullet
\save*++!U{1}\restore
\ar[l]
}
\hskip 30pt
\xymatrix@R=0.3in@C=0.4in{&&&{\bullet}\save*++!L{1}\ar[d]\restore\\
{\circ}
\save*++!U{0}\restore
\ar[r]
&
\circ
\save*++!U{0}\restore
\ar[r]
&
\circ
\save*++!U{0}\restore
\ar[r]
&
\bullet
\save*++!U{1}\restore
&
\bullet
\save*++!U{1}\restore
\ar[l]
&
\bullet
\save*++!U{1}\restore
\ar[l]
}
$$
$$
\xymatrix@R=0.3in@C=0.4in{&&&{\bullet}\save*++!L{1}\ar[d]\restore\\
{\circ}
\save*++!U{0}\restore
\ar[r]
&
\bullet
\save*++!U{1}\restore
\ar[r]
&
\bullet
\save*++!U{1}\restore
\ar[r]
&
\bullet
\save*++!U{1}\restore
&
\circ
\save*++!U{0}\restore
\ar[l]
&
\circ
\save*++!U{0}\restore
\ar[l]
}
\hskip 30pt
\xymatrix@R=0.3in@C=0.4in{&&&{\bullet}\save*++!L{1}\ar[d]\restore\\
{\circ}
\save*++!U{0}\restore
\ar[r]
&
\circ
\save*++!U{0}\restore
\ar[r]
&
\bullet
\save*++!U{1}\restore
\ar[r]
&
\bullet
\save*++!U{1}\restore
&
\bullet
\save*++!U{1}\restore
\ar[l]
&
\circ
\save*++!U{0}\restore
\ar[l]
}
$$
\vskip 10pt
The last root gives rise to the matrix
$$
\begin{array}{|c|c|c|c|}
\hline
C&-cI_4&0&0\\
\hline
0&-dI_4&D&0\\
\hline
0&-fI_4&0&F\\
\hline
\end{array}
$$
and assuming $c\neq 0$, column and row operations transform this into
$$
\begin{array}{|c|c|c|c|}
\hline
0&-cI_4&0&0\\
\hline
-\frac{d}{c}C&0&D&0\\
\hline
-\frac{f}{c}C&0&0&F\\
\hline
\end{array}
$$
If also $df\neq 0$, then this determinant vanishes if and only if that of
$$
\begin{array}{|c|c|c|}
\hline
-C&D&0\\
\hline
-C&0&F\\
\hline
\end{array}
$$
vanishes, which is the case when $\im\,C\cap\im\,D\cap\im\,F\neq 0$; this 
can be seen by noting that if $Cu=Dv=Fw$ then the vector 
$(u,v,w)^t$ lies in its kernel, and vice versa. \\

\ni The sixth and last component of $D$ is given by the vanishing of the 
semi-invariant arising from the root represented by $W$ in the diagram
$$
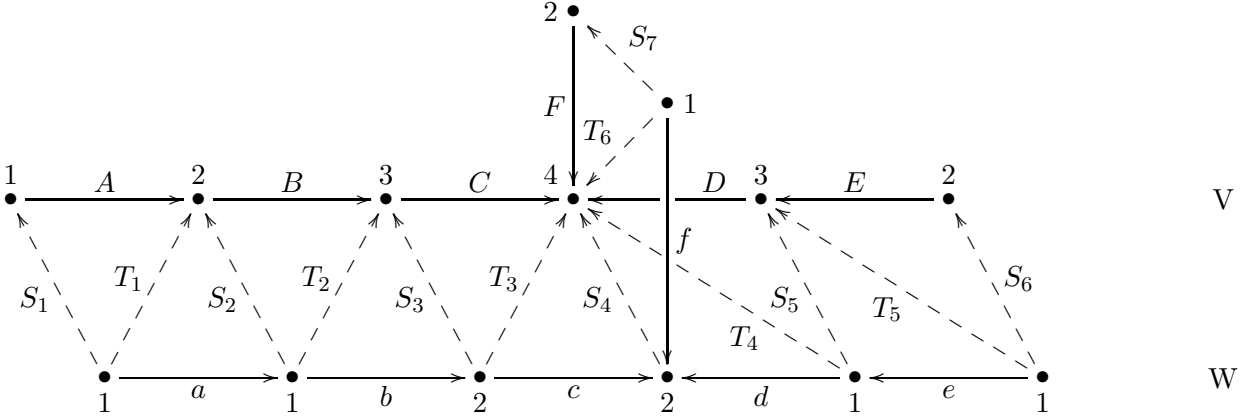
\xymatrix{
&&&&&&\bullet
\save*++!R{2}\restore
\ar[dd]_{\ds F}\\
&&&&&&&\bullet
\save*++!L{1}\restore
\ar@{-->}[ul]
_{\ds S_7}
\ar@{-->}[dl]
_{\ds T_6}
\ar[ddd]^{\ds f}
\\
\bullet
\save*++!D{1}\restore
\ar[rr]^{\ds A}
&&
\bullet
\save*++!D{2}\restore
\ar[rr]^{\ds B}
&&
\bullet
\save*++!D{3}\restore
\ar[rr]^{\ds C}
&&
\bullet
\save*++!DR{4}\restore
&&
\bullet
\save*++!D{3}\restore
\ar[ll]|\hole_<<<<<{\ds D}
&&
\bullet
\save*++!D{2}\restore
\ar[ll]_{\ds E}
&&&\txt{V}
\\\\
&
\bullet
\save*++!U{1}\restore
\ar@{-->}[uul]^{\ds S_1}
\ar@{-->}[uur]^{\ds T_1}
\ar[rr]_{\ds a}
&&
\bullet
\save*++!U{1}\restore
\ar@{-->}[uul]^{\ds S_2}
\ar@{-->}[uur]^{\ds T_2}
\ar[rr]_{\ds b}
&&
\bullet
\save*++!U{2}\restore
\ar@{-->}[uul]^{\ds S_3}
\ar@{-->}[uur]^{\ds T_3}
\ar[rr]_{\ds c}
&&
\bullet
\save*++!U{2}\restore
\ar@{-->}[uul]^{\ds S_4}
&&
\bullet
\save*++!U{1}\restore
\ar[ll]^{\ds d}
\ar@{-->}[llluu]|>>>>>>>>>>>>>\hole^<<<<<<<<<<<<<{\ds T_4}
\ar@{-->}[uul]^{\ds S_5}
&&
\bullet
\save*++!U{1}\restore
\ar@{-->}[llluu]^{\ds T_5}
\ar@{-->}[luu]_{\ds S_6}
\ar[ll]^{\ds e}
&&\txt{W}
}
$$

\ni The resulting determinant is
\vskip 10pt
$$
\begin{array}{|c|c|cc|cc|c|c|c|}
\hline
A  &  -aI_2  &  0  &  0  &  0  &  0  &  0 & 0 & 0\\
\hline
0  &  B  & -b_{11}I_3  &-b_{21}I_3    &  0  &  0  &  0 & 0 & 0\\
\hline
0&0&C&0&-c_{11}I_4&-c_{21}I_4&0&0&0\\
0&0&0&C&-c_{12}I_4&-c_{22}I_4&0&0&0\\
\hline
0&0&0&0&-d_{11}I_4&-d_{21}I_4&D&0&0\\
\hline
0&0&0&0&0&0&-eI_3&E&0\\
\hline
0&0&0&0&-f_{11}I_4&-f_{21}I_4&0&0&F\\
\hline
\end{array}
$$
\vskip 10pt
\ni where the columns and rows refer, in this order, to the maps $S_1\ld,S_7$ 
and $T_1,\ld,T_6$ respectively.
Row and column operations, and the deletion of rows and columns containing
only an invertible matrix, transform this to the matrix 
\vskip 10pt
$$
\begin{array}{|c|c|c|c|c|c|}
\hline
\frac{1}{b_{11}a}CBA & 0&(c_{12}b_{21}-c_{11}b_{11})I_4&(c_{22}b_{21}-
c_{21}b_{11})I_4&0&0\\
\hline
0&C&-c_{12}I_4&-c_{22}I_4&0&0\\
\hline
0&0&-d_{11}I_4&-d_{21}I_4&\frac{1}{e}DE&0\\
\hline
0 &  0 & -f_{11}I_4&-f_{21}I_4& 0 & F\\
\hline
\end{array}
$$
\vskip 10pt
\ni and now permuting columns brings it to the form
\vskip 10pt
$$
\begin{array}{|c|c|c|c|c|c|}
\hline
CBA&0&0&0&\l_1I_4&\mu_1I_4\\
\hline
0&C&0&0&\l_2I_4&\mu_2I_4\\
\hline
0&0&DE&0&\l_3I_4&\mu_3I_4\\
\hline
0&0&0&F&\lambda_4I_4&\mu_4I_4\\
\hline
\end{array}
$$
\vskip 10pt
\ni where the $\l_i$ and $\mu_j$ are polynomials in the coefficients 
$a,b,\ld$ of the representation $W$, and we have multiplied some rows
and columns by other such polynomials to simplify the expression
(since we choose a generic $W$ in $\Rep(Q,\e)$ to obtain
the polynomial $C^W$, this multiplication has the effect only of 
multiplying $C^W$ by a scalar).\\

\ni The geometrical significance of the vanishing of the determinant is 
that the three lines $\im\,DE\cap\im\,C,\ \im\,F\cap\im\,C$ and $\im\,CBA$
fail to span $\im\,C$. The reader will note the similarity in the geometric 
description of the last two semi-invariant factors. This can be understood 
by looking at their weights. They are given
by
\begin{align*}
\begin{matrix}
&&&-1\\
0&0&-1&2&-1&0
\end{matrix}\qquad\text{and}\qquad
\begin{matrix}
&&&-1\\
-1&0&-1&2&0&-1
\end{matrix}
\end{align*}
According to Derksen and Weyman \cite[p.477, Step 2]{derksenweyman}, 
if the weight of $W$ is not sincere, as in these cases, one may simplify 
the calculation by removing successively nodes not in the support, 
adding instead one arrow for each pair of ingoing and outgoing 
arrows. In the first case at hand, this produces a weight with support 
a Dynkin quiver of type $D_{4}$, in the second a weight of type $D_{5}$. 
For the first four orthogonal roots listed, the type of the weight equals 
$A_{3}$, explaining the similarity in the description of the corresponding 
semi-invariants. 
Once one has modified the quiver in this fashion, one can then simply 
calculate the corresponding semi-invariant on the new quiver, where 
one drops from $\d$ as well the nodes not in the support of the weight, 
and substituting at the end the actual composition of the maps along 
each pair of ingoing and outgoing arrow into the resulting semi-invariant. 
Revisiting, for example, the first orthogonal root considered above and 
its corresponding weight of type $A_{3}$; see e.g. the table below; 
it becomes thus transparent that the semi-invariant obtained, 
$\det[CBA|D]$, has indeed to be a polynomial in the entries of 
$CBA$ and $D$.

We can summarize the information gathered so far for the discriminant 
in the representation variety of the highest root of $E_{7}$ in the 
given orientation through the following table, where we list the opposite 
of the weights to display fewer minus signs:
$$
\begin{array}{|c|c|c|c|c|}
\hline
{\text{\sc Polynomial}}&\text{\sc Deg}&\text{\sc Root$^{\perp {\bf d}}$}&
\text{\sc $-$Weight}&\begin{matrix}\text{\sc Type}\\
\text{(Root, Weight)}
\end{matrix}\\
\hline
P_{1}=\det[CBA | D]&6 
&
\begin{matrix}
&&&0\\
1&1&1&1&1&0
\end{matrix}
&
\begin{matrix}
&&&0\\
1&0&0&-1&1&0
\end{matrix}
& (A_{5},A_{3})\\
\hline
P_{2}=\det[CB | DE]&8 
&
\begin{matrix}
&&&0\\
0&1&1&1&1&1
\end{matrix}
&
\begin{matrix}
&&&0\\
0&1&0&-1&0&1
\end{matrix}
& (A_{5},A_{3})\\
\hline
P_{3}=\det[F | DE]&6 
&
\begin{matrix}
&&&1\\
0&0&0&1&1&1
\end{matrix}
&
\begin{matrix}
&&&1\\
0&0&0&-1&0&1
\end{matrix}
& (A_{4},A_{3})\\
\hline
P_{4}=\det[CB | F]&6 
&
\begin{matrix}
&&&1\\
0&1&1&1&0&0
\end{matrix}
&
\begin{matrix}
&&&1\\
0&1&0&-1&0&0
\end{matrix}
& (A_{4},A_{3})\\
\hline
P_{5}=\det\
\left[
\begin{matrix}
-C&D&0\\
-C&0&F
\end{matrix}
\right]
&8
&
\begin{matrix}
&&&1\\
0&0&1&1&1&0
\end{matrix}
&
\begin{matrix}
&&&1\\
0&0&1&-2&1&0
\end{matrix}
& (D_{4},D_{4})\\
\hline
P_{6}&12
&
\begin{matrix}
&&&1\\
1&1&2&2&1&1
\end{matrix}
&
\begin{matrix}
&&&1\\
1&0&1&-2&0&1
\end{matrix}
& (E_{7},D_{5})\\
\hline
\Delta = (\text{unit})P_{1}\cdots P_{6}&46 
&
&
\begin{matrix}
&&&4\\
2&2&2&-8&2&3
\end{matrix}&
\\
\hline
\end{array}
$$
}
\end{example}
The following interlude will allow us to find the equations for 
semi-invariants such as $P_{5}$ or $P_{6}$ above in a more direct form, 
using some commutative algebra.
\section{An Interlude from Commutative Algebra}

Let $0\to M\xrightarrow{j}  R^{m+a}\xrightarrow{\pp} R^{a}\xrightarrow{p} 
T\to 0$ 
be an exact sequence of modules over a commutative normal (and noetherian) 
domain $R$, with integers $m,a >0$, and $T$ a torsion $R$--module. 
Assume given moreover an $R$--linear map $\psi:R^{m+a}\to R^{m}$. 
The module $M$ has a (constant) rank, equal to $m$, and its 
{\em determinant\/} is by definition the reflexive $R$--module $\det M = 
(\Lambda_{R}^{m}M)^{\vee\vee}$, where $(-)^{\vee}$ denotes 
the $R$--dual module. In words, $\det M$ is the reflexive hull of
 the $m^{th}$ exterior power of $M$ over $R$. It is isomorphic to $R$, and
the composition  $\psi j$ induces an $R$--linear map $\det(\psi j):R\cong 
\det M \to 
 \det R^{m}\cong R$. At issue now is to find a closed form for that determinant.
 
\begin{lemma} 
\label{lem:det}
The determinant of $\psi j$ satisfies $\det(\psi j) = \det(\psi, \pp)$.
\end{lemma}

\ni \proof
Consider the following commutative diagram whose rows are exact
 \begin{align*}
\xymatrix{
0\ar[r]&M\ar[r]^{j}\ar[d]_{\psi j} &R^{m+a}\ar[r]^{\pp}\ar[d]_{(\psi,\pp)}&
R^{a}\ar[r]^{p} \ar@{=}[d]&T\ar[r]\ar[d]&0\\
0\ar[r]&R^{m}\ar[r]^-{in_{1}}&R^{m}\oplus R^{a}\ar[r]^-{pr_{2}}&R^{a}\ar[r]&0
}
\end{align*}
The multiplicativity of the determinant shows first that $\det M\cong R$ and 
then yields $\det(\psi j) = \det(\psi,\pp)$.
\eop

\begin{example}
{\em We use this result to find a closed form for the semi-invariant 
$P_{6}$ for the highest weight of 
$E_{7}$ described in the last section. Namely, with the same notations 
as there, that invariant 
measures whether the three lines $\im\,DE\cap\im\,C,\ \im\,F\cap\im\,C$ 
and $\im\,CBA$ span 
$\im\, C$. To translate this into multilinear algebra, note that it is 
equivalent to say that the fibre 
product $X$ of $DE$ with $C$ over their common target, the fibre product 
$Y$ of $F$ with 
$C$ over the common target, and the image $Z$ of $BA$ do not span the 
domain of $C$. 
Each of $X,Y,Z$ is a rank one submodule of the domain of $C$, which is 
a free module 
of rank $3$ over $R$, the ring of the representation variety. We thus expect 
the corresponding 
invariant to be $\det[X|Y|Z]$, and the preceding lemma lets us make this 
precise:
In the following diagram, the top row is a direct sum of three short 
exact sequences of graded $R$--modules
\begin{align*}
\xymatrix{
0\ar[r]&
*{\begin{matrix}
X\\
\oplus\\
Y\\
\oplus\\
R(-2)
\end{matrix}}
\ar[dd]_{
{\begin{matrix}
(i_{1}, i_{2}, BA)
\end{matrix}}}
\ar[rrr]|{
{\begin{matrix}
(i_{1},  j_{1})\\
\oplus\\
(i_{2}, j_{2})\\
\oplus\\
\text{id}_{R}
\end{matrix}}}
&&&
*{\begin{matrix}
R^{3}\oplus R^{2}\\
\oplus\\
R^{3}\oplus R(-1)^{2}\\
\oplus\\
R(-2)
\end{matrix}}
\ar[dd]^{\displaystyle M}
\ar[rrr]|{
{\begin{matrix}
(C,-F)\\
\oplus\\
(C,-DE)\\
\oplus\\
0
\end{matrix}}}
&&&
*{\begin{matrix}
R(1)^{4}\\
\oplus\\
R(1)^{4}\\
\oplus\\
0
\end{matrix}}
\ar@{=}[dd]
\ar[r]&
0\\
\\
0\ar[r]&R^{3}\ar[rrr]^-{\displaystyle
in_{1}}&&&
*{\begin{matrix}
R^{3}\\
\oplus\\
R(1)^{4}\\
\oplus\\
R(1)^{4}
\end{matrix}}
\ar[rrr]^-{\displaystyle pr_{23}}&&&
*{\begin{matrix}
R(1)^{4}\\
\oplus\\
R(1)^{4}\\
\end{matrix}}\ar[r]&
0
}
\end{align*}
where the maps $i_{1}, i_{2}, j_{1}, j_{2}, in_{1}$ are the natural 
inclusions, 
$pr_{23}$ the projection onto the sum of second and third factor, and the 
matrix $M$ is of the form:
$$
\begin{array}{|c|c|c|c|c|c|}
\hline &R^{3}&R^{2}&R^{3}&R(-1)^{2}&R(-2)\\
\hline
R^{3}&I&0&I&0&BA\\
\hline
R(1)^{4}&C&-F&0&0&0\\
\hline
R(1)^{4}&0&0&C&-DE&0\\
\hline
\end{array}
$$
The desired semi-invariant is now $\det(i_{1},i_{2},BA)$, which equals 
the determinant of $M$ in view of the lemma above. Subtracting 
(a multiple of) the first column from the third and fifth results in the 
following simpler form
$$
\begin{array}{|c|c|c|c|c|c|}
\hline &R^{3}&R^{2}&R^{3}&R(-1)^{2}&R(-2)\\
\hline
R^{3}&I&0&0&0&0\\
\hline
R(1)^{4}&C&-F&-C&0&-CBA\\
\hline
R(1)^{4}&0&0&C&-DE&0\\
\hline
\end{array}
$$
whence the desired semi-invariant is seen to be the determinant of 
an $8\times 8$ matrix,
\begin{align*}
P_{6} = \det\left[
\begin{matrix}
F&C&0&CBA\\
0&C&-DE&0
\end{matrix}
\right]
\end{align*}
whose degree can be read off to be $12$ as stated in the table above.
}
\end{example}
\section{The case of $E_{8}$ with the centre as only sink}\label{e8}
As our final example, we determine the discriminant in the representation 
variety that belongs to the highest root of the Dynkin quiver of type $E_{8}$
with all arrows oriented towards the central 
trivalent vertex: 
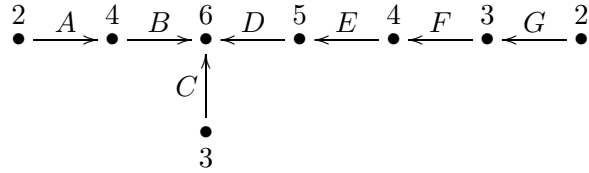
\[
\xymatrix{
\bullet\save*++!D{ 2} \restore
\ar[r]^{\displaystyle A} &
\bullet \save*++!D{ 4} \restore
\ar[r]^{\displaystyle B} &
\bullet \save*++!D{ 6} \restore &
\bullet \save*++!D{5} \restore
\ar[l]_{\displaystyle D} &
\bullet \save*++!D{ 4} \restore
\ar[l]_{\displaystyle E} &
\bullet \save*++!D{3} \restore
\ar[l]_{\displaystyle F} &
\bullet \save*++!D{2} \restore
\ar[l]_{\displaystyle G} \\
&&\bullet \save*++!U{ 3} \restore
\ar[u]^{\displaystyle C}
}\]
The capital letters $A,...,G$ stand for the corresponding matrices of 
independent indeterminates, and the coordinate ring of $\Rep(E_{8}.{\bf d})$ 
is $R=K[A,B,C,D,E,F,G]$, a polynomial ring in $118=120-2$ variables, 
where $120$ is the order of the binary icosahedral group.

We will also need below three additional auxiliary vertices, denoted by 
${\circ}$, and 
corresponding maps $X,Y,Z$, as indicated by the dashed arrows here:
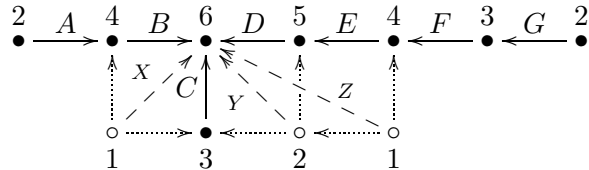
\[
\xymatrix{
\bullet\save*++!D{ 2} \restore
\ar[r]^{\displaystyle A} &
\bullet \save*++!D{ 4} \restore
\ar[r]^{\displaystyle B} &
\bullet \save*++!D{ 6} \restore &
\bullet \save*++!D{5} \restore
\ar[l]_{\displaystyle D} &
\bullet \save*++!D{ 4} \restore
\ar[l]_{\displaystyle E} &
\bullet \save*++!D{3} \restore
\ar[l]_{\displaystyle F} &
\bullet \save*++!D{2} \restore
\ar[l]_{\displaystyle G} \\
&\circ \save*++!U{1} \restore
\ar@{.>}[u]\ar@{.>}[r]\ar@{-->}[ru]^{X}
&\bullet \save*++!U{ 3} \restore
\ar[u]^{\displaystyle C}
&\circ \save*++!U{2} \restore
\ar@{.>}[u]|!{[r];[lu]}\hole\ar@{.>}[l]\ar@{-->}[ul]^{Y}
&\circ \save*++!U{1} \restore
\ar@{.>}[u]\ar@{.>}[l]\ar@{-->}[llu]_(.3){Z}
}\]
The map $X$ is the natural one from the fibre product of $B$ and $C$ to the 
central node. The fibre product itself is an $R$--module of rank $1$. 
The map $Y$ indicated above is the natural one from the fibre product 
of $D$ and $C$ to the central node. This fibre product has rank 2. 
Finally, the map $Z$ is the natural one from the fibre product of $C$ and 
$DE$ to the central node. Again, the fibre product has rank 1.

The discriminant $D$ in question is of degree 118 and has $7$ irreducible 
components, thus, its canonical equation $\Delta$ is a product of $7$ 
irreducible polynomials 
$P_{i}$ in the entries of the $7$ matrices $A$ through $G$. 
Moreover, we obtain from \ref{lem:weights} that it is a 
semi-invariant belonging to the weight
\begin{align*}
\begin{matrix}
-4&-4&12&-2&-2&-2&-3\\
&&-6
\end{matrix}
\end{align*}

We can spot immediately three semi-invariants:
\begin{align*}
P_{1}= \det[BA|DE]\,, P_{2}= \det[C|DEF]\,, P_{3}= \det[B|DEFG]\,,
\end{align*}
each of degree 12 and belonging to weights of type $A_{3}$. 
The remaining four can be described thus
\begin{itemize}
\item The failure of $\im\, X = \im\,B\cap\im\, C$ and $\im\, D$ to 
generate the vector space at the central node.
According to \ref{lem:det}, the corresponding polynomial is 
$P_{4} = \det[X|D]$, the determinant of
\begin{align*}
\begin{array}{|c|c|c|c|}
\hline
&R(-1)^{4}&R(-1)^{3}& R(-1)^{5}\\
\hline
R^{6}&B&0&D\\
\hline
R^{6}&B&-C&0\\
\hline
\end{array}
\end{align*}
It is of degree 12. 
\item
The failure of $\im\, BA, \im\, X, \im\, DEF$ to generate the 
vector space at the central node.
Again using \ref{lem:det}, the corresponding polynomial is 
$P_{5}=\det[BA|X|DEF]$, the determinant of
\begin{align*}
\begin{array}{|c|c|c|c|c|}
\hline
&R(-2)^{2}&R(-1)^{4}&R(-1)^{3}& R(-3)^{3}\\
\hline
R^{6}&BA&B&0&DEF\\
\hline
R^{6}&0&B&-C&0\\
\hline
\end{array}
\end{align*}
It is of degree 20.
\item The failure of $\im\, BA, \im\, Y, \im\, DEFG$ to generate the 
vector space at the central node. The corresponding polynomial is 
$P_{6}=\det[BA|Y|DEFG]$, the determinant of
\begin{align*}
\begin{array}{|c|c|c|c|c|}
\hline
&R(-2)^{2}&R(-1)^{3}&R(-1)^{5}& R(-4)^{2}\\
\hline
R^{6}&BA&C&0&DEFG\\
\hline
R^{6}&0&C&-D&0\\
\hline
\end{array}
\end{align*}
It is also of degree 20. The three semi-invariants $P_{4}$ through $P_{6}$ 
belong to weights of type $D$, as can easily be seen from the geometric 
description. Now we turn to the last and biggest one:
\item The rank of $BA$ and $DEFG$ is $2$, that of $X, Z$ is $1$. Their 
images in the central vector space are thus expected to generate. The 
failure will be measured by the polynomial
$P_{7} =\det[BA|X|Z|DEFG]$, which is the determinant of 
\begin{align*}
\begin{array}{|c|c|c|c|c|c|c|}
\hline
&R(-2)^{2}&R(-1)^{4}&R(-1)^{3}& R(-1)^{3}&R(-2)^{4}&R(-4)^{2}\\
\hline
R^{6}&BA&C&0&C&0&DEFG\\
\hline
R^{6}&0&C&-D&0&0&0\\
\hline
R^{6}&0&0&0&C&-DE&0\\
\hline
\end{array}
\end{align*}
It is of degree 30 and its weight is of type $E_{6}$.
\end{itemize}
We summarize the results again in a table:

$$
\begin{array}{|c|c|c|c|c|}
\hline
{\text{\sc Polynomial}}&\text{\sc Deg}&\text{\sc Root$^{\perp {\bf d}}$}&
\text{\sc $-$Weight}&\begin{matrix}\text{\sc Type}\\
\text{(Root, Weight)}
\end{matrix}\\
\hline
P_{1}= \det\left(BA | DE\right)&12 
&
\begin{matrix}
1&1&1&1&1&0&0\\
&&0
\end{matrix}
&
\begin{matrix}
1&0&-1&0&1&0&0\\
&&0
\end{matrix}
& (A_{5},A_{3})\\
\hline
P_{2}= \det\left(C | DEF\right)&12 
&
\begin{matrix}
0&0&1&1&1&1&0\\
&&1
\end{matrix}
&
\begin{matrix}
0&0&-1&0&0&1&0\\
&&1
\end{matrix}
& (A_{5},A_{3})\\
\hline
P_{3}= \det\left(B | DEFG\right)&12 
&
\begin{matrix}
0&1&1&1&1&1&1\\
&&0
\end{matrix}
&
\begin{matrix}
0&1&-1&0&0&0&1\\
&&0
\end{matrix}
& (A_{6},A_{3})\\

\hline
P_{4}= \det\left( X | D\right)&12 
&
\begin{matrix}
0&1&1&1&0&0&0\\
&&1
\end{matrix}
&
\begin{matrix}
0&1&-2&1&0&0&0\\
&&1
\end{matrix}
& (D_{4},D_{4})\\
\hline
P_{5}= \det\left(BA | X | DEF\right)&20 
&
\begin{matrix}
1&2&2&1&1&1&0\\
&&1
\end{matrix}
&
\begin{matrix}
1&1&-2&0&0&1&0\\
&&1
\end{matrix}
& (E_{7},D_{5})\\
\hline
P_{6}= \det\left(BA | Y | DEFG\right)&20 
&
\begin{matrix}
1&1&2&2&1&1&1\\
&&1
\end{matrix}
&
\begin{matrix}
1&0&-2&1&0&0&1\\
&&1
\end{matrix}
& (E_{8},D_{5})\\
\hline
P_{7}= \det\left(BA | X|Z | DEFG\right)&30 
&
\begin{matrix}
1&2&3&2&2&1&1\\
&&2
\end{matrix}
&
\begin{matrix}
1&1&-3&0&1&0&1\\
&&2
\end{matrix}
& (E_{8},E_{6})\\
\hline
\Delta = (\text{unit})P_{1}\cdots P_{7}&118 
&
&
\begin{matrix}
4&4&-12&2&2&2&3\\
&&6
\end{matrix}&
\\
\hline
\end{array}
$$

\begin{theorem}
The above table is correct.
\end{theorem}

\begin{proof}
Inspection shows that each of the polynomials $P_{i}$ is a semi-invariant and that its weight is as listed. The indicated weights are easily seen to be linearly independent, and add up to the weight of $\Delta$. Thus, their product must describe the discriminant up to a unit.
\end{proof}


\begin{thebibliography}{cc}
\bibitem{ARS} M.Auslander, I.Reiten, and S.O.Smal$\o$,  Representation theory of Artin algebras. Corrected reprint of the 1995 original.
Cambridge Studies in Advanced Mathematics, 36. Cambridge University Press, Cambridge,  1997. xiv+425 pp.
\bibitem{bgp} I.N.Bernstein, I.M.Gel'fand and V.A.Ponomarev, Coxeter
functors and Gabriel's Theorem, Uspekhi Mat. Nauk. 28 (1973), 19-33,
translated in Russian Math. Surveys 28 (1973) 17-32
\bibitem{B} N.Bourbaki, Groupes et Alg\`ebres de Lie, Chapitres 4,5 et 6,
Masson, 1981
\bibitem{buchweitz} R.-O.Buchweitz, {\it in preparation}
\bibitem{buchweitzebeling} R.-O.Buchweitz and W.Ebeling, {\it in preparation}
\bibitem{damon1} J.N.Damon, Non-linear sections of nonisolated complete
intersections, in {\it New developments in singularity theory (Cambridge 2000)},
NATO Sci. Ser.II Math.Phys.Chem. 21, 405-445, Kluwer Acad. Publ., Dordrecht,
2001
\bibitem{damonleg} J.N.Damon, On the legacy of free divisors II: 
Free* divisors and complete intersections, Moscow Math.J. (3) no. 2, 
(2003) 361-395
\bibitem{derksenweyman} H.Derksen and J.Weyman, Semi-invariants of quivers and 
saturation for Littlewood-Richardson coefficients, J.Amer. Math.Soc. 13
(2000) 467-479
\bibitem{df} P.Donovan and M.R.Freislich, {\it The representation theory of
finite graphs and associated algebras}, Carleton Math. Lecture Notes No. 5, 
1973
\bibitem{gabriel} P.Gabriel, Unzerlegbare Darstellungen I, Manuscripta Math.
6 (1972) 71-103
\bibitem{gabriel1} P.Gabriel, R\'epresentation ind\'ecomposables, S\'eminaire
Bourbaki, Expos\'e 444 (1974), in Springer Lecture Notes 431 (1975) 143-169
\bibitem{cg} C.G.Gibson, {\it Singular points of smooth mappings}, Res. Notes
in Math. 25, Pitman, London, 1979
\bibitem{goryunov} V.Goryunov, Functions on space-curves, J.London Math.Soc
61 (2000), 807-822
\bibitem{gmns} M.Granger, D.Mond, A.Nieto, and M.Schulze, Linear free divisors,
{\it in preparation}
\bibitem{HU} D.Happel, S.Hartlieb, O.Kerner, and L.Unger, On perpendicular categories of stones over quiver algebras.
 Comment. Math. Helv.  71  (1996),  no. 3, 463--474.
\bibitem{kac} V.Kac, Infinite root systems, representations of graphs and
invariant theory, Invent. Math. 56 (1980), 57-92
\bibitem{kac2} V.Kac, Infinite root systems, representations of graphs and
invariant theory II, J.Algebra 78 (1982), 57-92 
\bibitem{krre}H.Kraft and Ch.Riedtmann, Geometry of representations of
quivers, in {\it Representations of Algebras, Durham 1985}, P.Webb (Ed.),
London Math. Soc. Lecture Notes 116, 1986, 109-146
\bibitem{looijenga} E.J.N.Looijenga, {\it Isolated singular points on complete
intersections}, London Math.Soc. Lecture Notes in Math. 77, Cambridge
University Press, 1984
\bibitem{matherIV}J.N.\,Mather, Stability of $C^{\infty}$ maps IV: Classification
of stable maps by $R$-algebras, Publ.Math. I.H.E.S. 37(1969) 223-248
\bibitem{ms} D.Mond and D.van Straten, Milnor number equals Tjurina number
for functions on space curves, J.London Math. Soc. (2) 63 (2001), 177-187
\bibitem{mwik} D.Mond and R.Wik, Not all codimension 1 germs have good real
pictures, in {\it Real and complex singularities}, D.Mond and M.Saia (eds.),
Marcel Dekker 2003, 189-200
\bibitem{naz} L.A.Nazarova, Representations of quivers of infinite type,
Math. USSR Izvestija, Ser. Mat. 37 (1973) 752-791
\bibitem{orlikterao} P.Orlik and H.Terao, {\it Arrangements of hyperplanes}, 
Grundlehren der Math. Wissenschaft 300, Springer Verlag, 1992
\bibitem{ringelalg} C.M.Ringel, Representations of K-species and bimodules,
J.Alg. 41 (1976) 269-302
\bibitem{ringelinv} C.M.Ringel, Rational invariants of tame quivers,
Invent. Math. 58 (1980) 217-239
\bibitem{saito} K.Saito, Theory of logarithmic differential forms and 
logarithmic vector fields, J.Fac.Sci. Univ.Tokyo Sect. Math. 27 (1980),
265-291
\bibitem{saki} M.Sato and T.Kimura, A classification of irreducible 
prehomogeneous vector spaces and their relative invariants, Nagoya J.Math.
65 (1977) 1-155
\bibitem{schofield} A.Schofield, Semi-invariants of quivers, J.London
Math. Soc. (2) 43 (1991) 385-395
\bibitem{schofield2} A.Schofield, General representations of quivers, Proc.
London Math. Soc. (3) 65 (1992) 46-64
\bibitem{duco} D.van Straten, A note on the discriminant of a space curve,
Manucripta Math 87 (1995), 167-177
\end{thebibliography}




{\sc Dept.\ of Math., University of
Tor\-onto, Tor\-onto, Ont.\ M5S 2E4, Canada}\\

{\it E-mail address:}\ \tt{ragnar@math.utoronto.ca}\\

{\sc Mathematics Institute, University of Warwick, Coventry CV4 7AL,
England}\\

{\it E-mail address:}\ \tt{mond@maths.warwick.ac.uk}
\end{document}